\documentclass{article}
\usepackage[version=4]{mhchem}
\usepackage[english]{babel}
\usepackage[letterpaper,top=2cm,bottom=2cm,left=3cm,right=3cm,marginparwidth=1.75cm]{geometry}
\usepackage{diagbox}
\usepackage{amsmath}
\usepackage{appendix}
\usepackage{ulem}
\usepackage[nobysame]{amsrefs}
\usepackage{amssymb,xcolor}
\usepackage{mathrsfs}
\usepackage{graphicx}
\usepackage{subfig}
\usepackage{float}
\usepackage{epsf}
\usepackage{color}
\usepackage{mathtools}
\usepackage[linesnumbered,ruled]{algorithm2e}
\usepackage{multicol}
\usepackage{multirow}
\usepackage[font=scriptsize]{caption}
\usepackage{url}
\usepackage{cases}
\usepackage{mhchem}
\usepackage{amsthm}

\newtheorem{thm}{Theorem}[section]

\newtheorem{rmk}{Remark}[section]

\title{Crossover from ballistic transport to normal diffusion: a kinetic view}
\author{Zhe Xue
\thanks{School of Mathematical Sciences, Institute of Natural Sciences, MOE-LSC, Shanghai Jiao Tong University, Shanghai, 200240, P.R. China. xuezhe\_017@sjtu.edu.cn}
\and
Weiran Sun
\thanks{Department of Mathematics, Simon Fraser University, Burnaby BC V5A 1S6, Canada. weirans@sfu.ca}
\and
Zhennan Zhou
\thanks{Institute for Theoretical Sciences, Westlake University, Hangzhou, 310030, P. R. China. zhouzhennan@westlake.edu.cn} \and 
Min Tang
\thanks{School of Mathematical Sciences, Institute of Natural Sciences, MOE-LSC, CMA-Shanghai, Shanghai Jiao Tong University, Shanghai, 200240, P.R. China. tangmin@sjtu.edu.cn}}
\date{}

\begin{document}
\maketitle

\begin{abstract}
     The crossover between dispersion patterns has been frequently observed in various systems. Inspired by the pathway-based kinetic model for \textit{E. coli} chemotaxis that accounts for the intracellular adaptation process and noise, we propose a kinetic model that can exhibit a crossover from ballistic transport to normal diffusion at the population level. 
     At the particle level, this framework aligns with a stochastic individual-based model. Using numerical simulations and rigorous asymptotic analysis, we demonstrate this crossover both analytically and computationally. Notably, under suitable scaling, the model reveals two distinct limits in which the macroscopic densities exhibit either ballistic transport or normal diffusion.


     
\end{abstract}

\section{Introduction}

The random movement of a particle was described by Brown in his foundational studies of pollen grains and other microscopic particles \cite{Brownian1828}. Since then, various tools have been developed to model and measure how particles disperse over time. A standard statistical quantity used to quantify this dispersion is the expected or mean squared displacement (MSD) at each time. 
Since particles change their positions through stochastic processes, MSD serves as a metric to distinguish between normal and anomalous diffusion. 
More precisely, let $x(t)$ be the displacement of a particle at time $t$. In cases where
\begin{equation*}
{\rm MSD} = \left \langle x^2(t) \right \rangle \sim K_{\alpha} t^{\alpha},
\end{equation*}
the particles undergo a normal diffusion 
if the exponent $\alpha$ is equal to 1. 
When the scaling exponent $\alpha$ deviates from 1, the motion is considered subdiffusive for $0 < \alpha < 1$, superdiffusive for $1 < \alpha < 2$, and ballistic when $\alpha = 2$. 
%
Each of these modes of dispersion has been observed in various physical and biological systems. For instance, as discussed in \cite{di2021subdiffusive}, densely packed and heterogeneous structures result in subdiffusion for lipids and proteins; superdiffusion in several cellular systems has been documented, which is attributed to active motion \cites{chen2015memoryless, Song2018NeuronalMR}. In semiconductor physics, the lack of scattering centers facilitates the ballistic transport of electrons or holes \cites{Shurballistic, Freyballistic}. For more examples, one may consult the review articles \cites{Höfling_2013, norregaard2017manipulation, metzler2016non}.

An interesting phenomenon of a crossover between different modes of dispersion has been frequently observed in complex systems \cites{PhysRevLett.120.248101, Molina-Garcia_2018, PhysRevLett.109.188103, PhysRevE.89.062126, Dentz2003TimeBO, PhysRevLett.108.230602}. 
 In particular, the transition from the ballistic transport to normal diffusion has been observed in various physical and biological systems. For example, under symmetric space-time reflection, wave transport can experience an abrupt shift from ballistic to diffusive behavior at a critical temporal juncture \cite{eichelkraut2013mobility}. Moreover, in scenarios characterized by weak interactions and finite energy density, fermion-fermion scattering leads to a crossover from early-time ballistic to late-time diffusive behavior \cite{crossoverballisticLloyd}. Another noteworthy example is animal migration studied in \cite{tilles2016animals}, where the model developed shows that the animals' movement trajectories evolve from an initial ballistic transport phase to a long-term normal diffusive phase. 


Several theoretical frameworks have been introduced to model the dispersion of particles. In \cites{MontrollWeiss, ScherMontroll}, Montroll and Weiss introduced the Continuous Time Random Walk (CTRW) model to capture the stochastic nature of diffusion processes. In the L\'{e}vy flight model described in \cite{levywalks}, particles make alternatively instantaneous jumps and waiting pauses; both the jump lengths and the durations of waiting events are random variables. When the distributions of these random variables transit from power-law to exponential decay, the particles, at the population level, exhibit a crossover from ballistic transport to normal diffusion \cite{LoopyLF}. Similarly, when random walkers with a constant speed alternately perform ballistic transport with a random duration time and reset movement directions randomly,  a crossover from ballistic transport to normal diffusion can also be observed by assigning an appropriate run length distribution \cite{BtoNLW}. However, all these crossovers are phenomenological and include only statistical information of particle trajectories.


Our goal in this paper is to provide an intracellular mechanism to generate the crossover from ballistic transport to normal diffusion in the physical space. 
More specifically, we focus on three aspects: 
 first, by simulating the random process of the internal state of the particle, we justify the particular choices of the tumbling frequency and adaptation function 
 (see Section~\ref{two-state}). 
Second, we set up a particle model incorporating the internal state dynamics and spatial movement of the cell, and numerically demonstrate a crossover from ballistic transport to normal diffusion across various adaptation times (see Section~\ref{sec:kinetic}).
Third, we formulate a kinetic equation from the particle model and analytically prove that as the adaptation time increases, the limits of the kinetic equation transition from self-similar ballistic transport to normal diffusion (see Sections~\ref{sec:results} and~\ref{sec:proof}).

This work is inspired by a state-dependent and spatially homogeneous jump model for a single particle \cite{XZZT}, where the probability of a particle jumping out of a given state is controlled by the intracellular signaling pathway. When intracellular adaptation is slow and noise is strong, the waiting time distribution transitions from a power-law to an exponential decay. 
The focus of the current paper is to include the spatially inhomogeneous dynamics of particle movement and to provide an intracellular mechanism for the crossover from ballistic transport to normal diffusion at the population level. 

Various dispersion limits, including normal and fractional diffusion, have been derived based on similar pathway-based kinetic models for Escherichia coli (\textit{E. coli}) chemotaxis ~\cites{Perthame2018, Xue2021, SiTangYang, PerthameSunTangShugo}. In these previous works, the tumbling frequency has to be strictly positive for normal diffusion~\cites{SiTangYang, PerthameSunTangShugo}, or degenerate at zero or infinity for a fractional diffusion limit~\cites{Perthame2018,Xue2021}. In the present work, we allow the tumbling frequency to degenerate in a finite interval. This enables particles to run for extended periods. As a consequence, a significant portion of the particles rarely tumble. The asymptotic behaviour in this case differs drastically from the earlier studies and the asymptotic analysis in the previous works does not apply, even at a formal level. Here we rely on using the explicit solution of the kinetic equation to derive the limiting solution and the equation it satisfies.

The structure of this paper is as follows. Section 2 constructs and simulates an individual-based model with intracellular signaling pathways. The observed crossover from ballistic transport to normal diffusion is verified through detailed statistical analysis based on numerical results. Section 3 introduces the kinetic model for the probability density function of the particles and presents the main analytical results. The proofs of these results are provided in Section 4, and we conclude with a discussion in Section 5. 

\section{The individual-based model}
In this section, we introduce an individual-based model that captures both the intracellular dynamics and the resulting population-level behavior. We begin by presenting an intracellular two-state model from \cite{tu2005white} for \textit{E. coli} chemotaxis. 
\textit{E. coli} cells alternately switch between run and tumble states, and the duration time of the run state has a transition from power-law decay to exponential decay, as pointed out in \cite{XZZT}.
By simulating the random process of the internal states of the particle, we identify necessary conditions for such transition to appear.
%
Furthermore, based on a simplified one-state model, we develop a kinetic PDE model that incorporates the run-and-tumble movements. This allows us to investigate the population-level behavior and how it is influenced by the internal adaptation dynamics. We will then conduct numerical simulations to explore the emerging dispersion patterns. It is important to note that the main purpose of this paper is to provide an intracellular mechanism for the transitional behavior of dispersion patterns, rather than to model a specific type of cells. Therefore, although the model is motivated by \textit{E. coli} cells, the specific forms of the functions and the values of the parameters may not exactly follow those of \textit{E. coli} cells.

\subsection{The two-state model and simplifications}\label{two-state}
\subsubsection{The two-state model}\label{two-statemodeltu}
\textit{E. coli} adeptly navigates its surroundings through run-and-tumble movements. During the run phase, \textit{E. coli} moves along a relatively straight line at a constant speed. In the tumble phase, the bacterium performs a random reorientation. The switching between the run and tumble phases of each \textit{E. coli} cell is modulated by the rotational directions of its flagella. Each \textit{E. coli} cell has roughly 6-8 flagella, each of which can rotate either clockwise (CW) or counterclockwise (CCW). The cell begins to run when more flagella rotate CCW, while it starts tumbling when more flagella rotate CW. In more realistic models, one must determine the relationship between the CCW and CW rotational phases of the 6-8 flagella and the cell's run-and-tumble movement phases \cites{MotileBehaviorofBacteria, TheRotaryMotorofBacterialFlagella, larsen1974change}. However, some models in the literature simplify the situation by assuming that each \textit{E. coli} cell has only one flagellum. When this flagellum rotates CCW, the cell runs; when it rotates CW, the cell tumbles.

In~\cite{tu2005white}, the authors proposed a two-state model to describe the switches between CCW and CW rotations. The switching from CCW to CW and CW to CCW are modeled by two stochastic Poisson processes, with rates denoted as $\Lambda_1$ and $\Lambda_2$, i.e.
\begin{equation*}
    \text{CCW} \ce{<=>[$~~~\Lambda_1~~~$][$~~~\Lambda_2~~~$]} \text{CW}.
\end{equation*}
 The switching rates are determined by the concentration of the response regulator protein CheY-P, which can bind to the flagellar motor and promote its rotation in the CCW direction. Let $Y_t$ be the concentration of CheY-P. As in \cite{tu2005white}, the switching rates from CCW to CW and from CW to CCW are given by 
\begin{equation}\label{TupaperLambda}
    \Lambda_1(Y_t) = t_0^{-1}\exp\left(\alpha_1 \frac{Y_t-\Bar{Y}}{\Bar{Y}}\right), \quad 
    \Lambda_2(Y_t) = t_1^{-1} \exp\left(\alpha_2 \frac{Y_t-\Bar{Y}}{\Bar{Y}}\right).
\end{equation}
Here $t_0$ and $t_1$ are the characteristic switching times; $\Bar{Y}$ is the mean value of $Y_t$; $\alpha_{1}$ and $\alpha_{2}$ are dimensionless constants. It is important to note that, in the model used in \cite{tu2005white}, $\alpha_1>0$, $\alpha_2<0$, and 
$\alpha_1$ is large. Therefore, according to the specific form in \eqref{TupaperLambda}, when $Y_t < \Bar{Y}$, $\Lambda_1 \ll \Lambda_2$ (and when $Y_t > \Bar{Y}$, $\Lambda_1 \gg \Lambda_2$), the flagella are more likely to switch from CW to CCW (and from CCW to CW). Moreover, $Y_t$ evolves according to an Ornstein-Uhlenbeck (OU) process 
\begin{equation}\label{TupaperY}
    dY_t = -\frac{Y_t-\Bar{Y}}{T_m} dt +\sigma dB_t = -\frac{\Delta Y}{T_m} dt +\sigma dB_t,
\end{equation}
where $\Delta Y = Y_t-\Bar{Y}$, $T_m$ denotes the CheY-P correlation time, $B_t$ is the standard Brownian motion and $\sigma$ is the intensity of the Brownian motion. 

We simulate the switching processes of the aforementioned two-state model, employing parameters similar to those in \cite{tu2005white}. In \cite{tu2005white}, the authors tested different noise strengths $\Delta_n$ defined as $\Delta_n \triangleq \alpha_1 \frac{\sigma \sqrt{T_m/2}}{\bar Y}$ and derive the corresponding CCW duration time distribution. 
In this paper, we use the same parameters $\alpha_1 = 10$, $\alpha_2 = -2$, $t_0 = 300$, $t_1 = 30$ and $T_m = 6000$ as in \cite{tu2005white}. For the noise strength, we select a large value of $\Delta_n = 50$ with $\bar{Y} = 5$ and $\sigma = 0.456$. In the simulation, we track $N = 100$ samples, with a total simulation time $6 \times 10^5$. Each sample includes information about its CheY-P concentration and the rotational state of the flagella (CCW or CW). Details of the numerical algorithm are described in Appendix \ref{twostateappendix}.
\begin{figure}[!h]  
    \centering
        \subfloat[]{\includegraphics
        [scale =0.5]{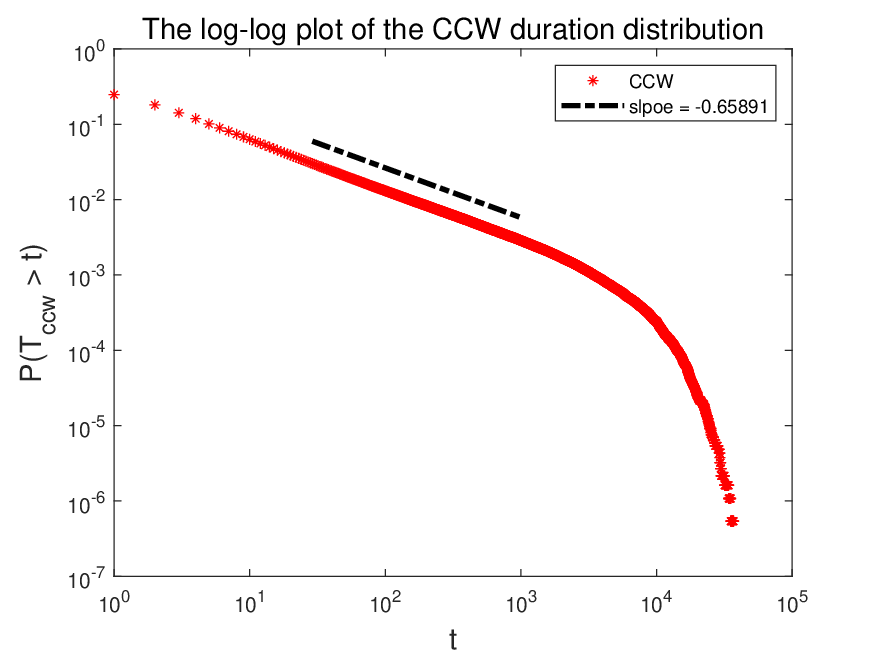}}
        \subfloat[]{\includegraphics
        [scale =0.5]{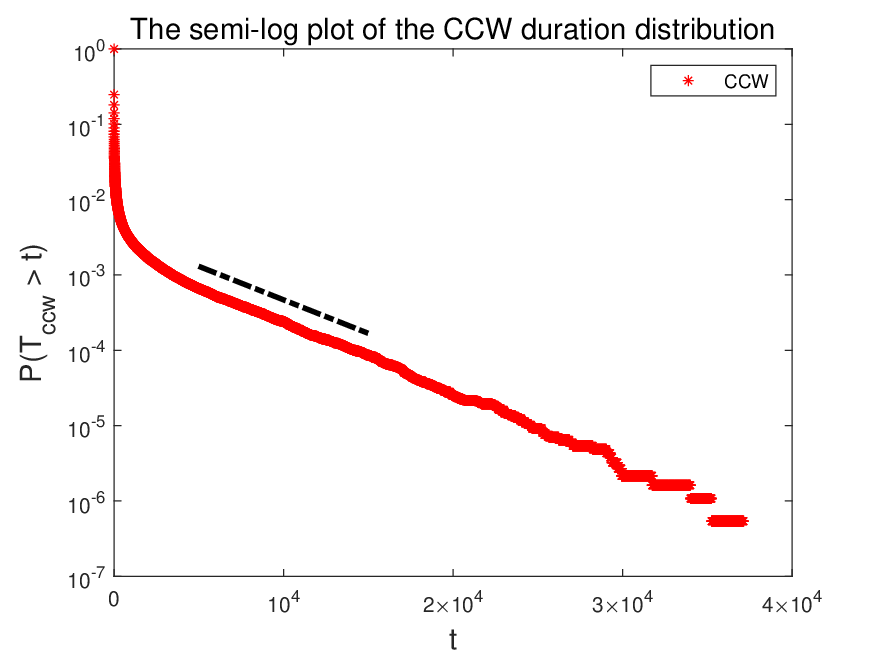}} 
        \caption{The distribution of the CCW duration time with $\alpha_1 = 10$, $\alpha_2 = -2$, $t_0 = 300$, $t_1 = 30$, $\bar{Y} = 5$, $T_m = 6000$ and $\sigma = 0.456$ in the two state model \eqref{TupaperLambda} and \eqref{TupaperY}. (a) The log-log plot of $P(T_{ccw} > t)$. The slope $-0.65891$ is fitted in the interval $[30, 10^3]$. (b) The semi-log plot $P(T_{ccw} > t)$. The straight dashed line is fitted in the interval $[5\times 10^3, 1.5\times10^4]$ which indicates the CCW duration time distribution decays exponentially.}
        \label{tupaperTccw}
\end{figure}
Since the flagella alternately rotate CCW or CW, the CCW (CW) duration time is the time period during which the flagella is in the CCW rotational state, which is a random variable. The CCW (CW) duration time is denoted by $T_{\text{ccw}}$ ($T_{\text{cw}}$). Figure \ref{tupaperTccw}(a) shows the probability distribution of $P(T_{\text{ccw}} > t)$. It can be seen that the CCW duration has a power-law decay in the intermediate interval $[30,10^3]$, as noted in~\cite{tu2005white}. We further observe in Figure \ref{tupaperTccw}(b) that the CCW duration time decays exponentially at the tail part $[3\times 10^3,10^4]$. Such transitional behavior in the CCW duration time distribution has only been studied recently in~\cite{XZZT}, both analytically and numerically. 

\subsubsection{The necessary ingredients for the transition in the CCW duration time}
\begin{figure}[!h]  
    \centering
        \subfloat[]{\includegraphics
        [scale =0.36]{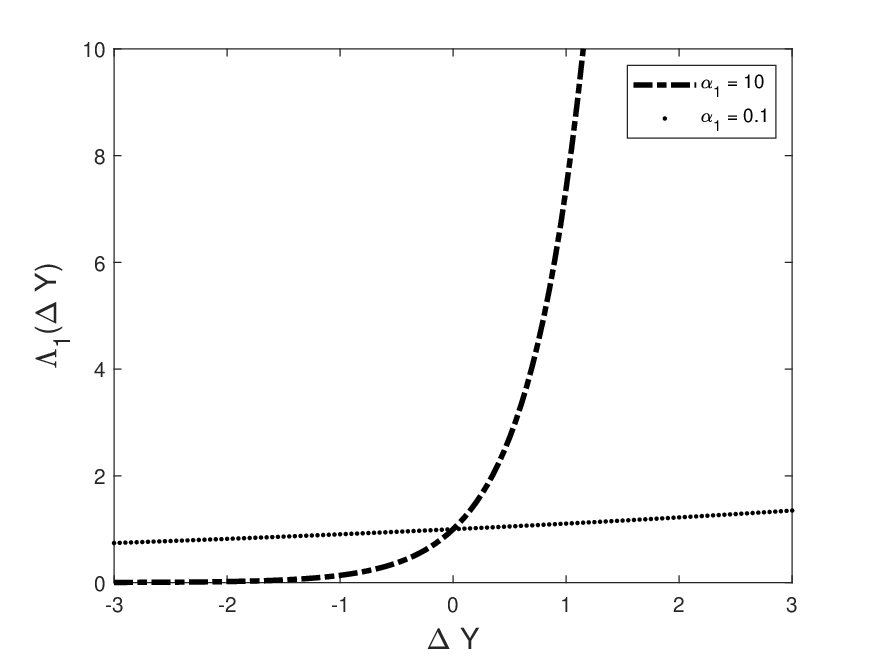}}
        \subfloat[]{\includegraphics
        [scale =0.36]{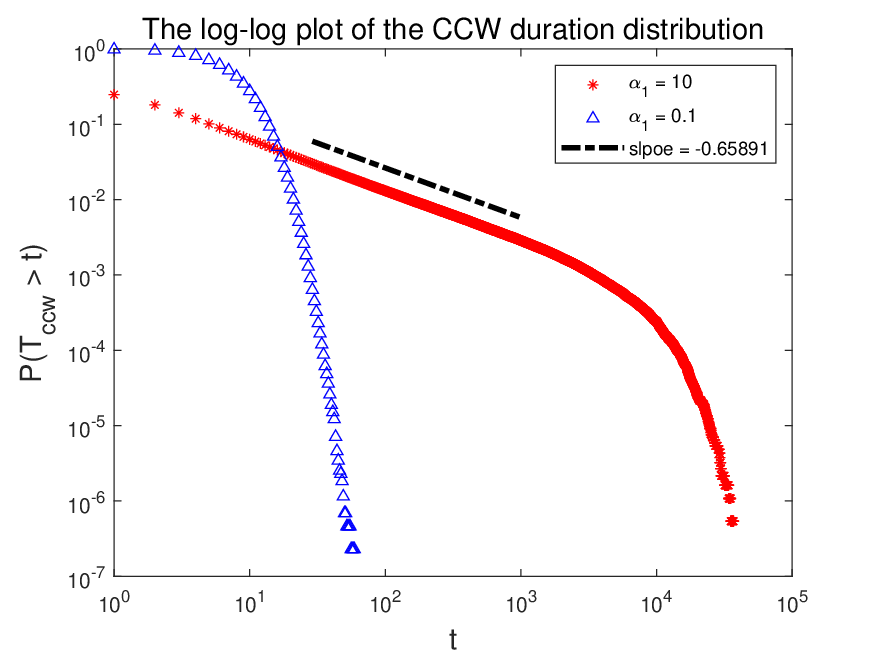}} 
        \subfloat[]{\includegraphics
        [scale =0.36]{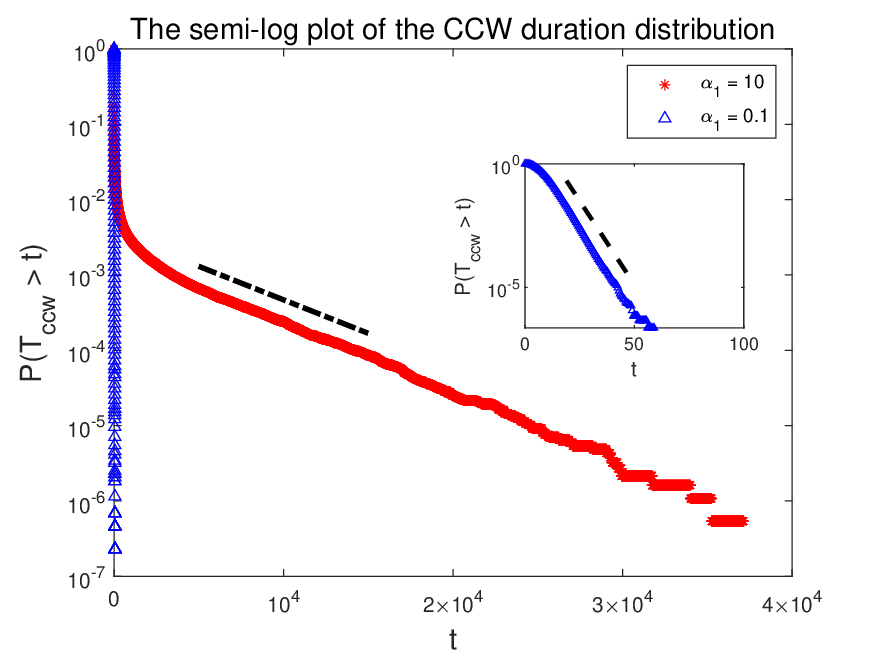}}
        \caption{ Comparison of different values of $\alpha_1$ in $\Lambda_1$ in \eqref{TupaperLambda}. (a) $\Lambda_1(\Delta Y) = \exp \left(\alpha_1 \frac{\Delta Y}{\bar{Y}}\right)$ for $\alpha_1 = 10$ (dash-dot line) and $\alpha_1 = 0.1$ (dotted line). (b) The log-log plot of $P(T_{ccw} > t)$ for $\alpha_1 = 10$ (red stars) and $\alpha_1 = 0.1$ (blue triangles); The slope $-0.65891$ is fitted over the interval $[30, 10^3]$. (c) The semi-log plot for $\alpha_1 = 10$ (red stars) and $\alpha_1 = 0.1$ (blue triangles). 
        The inset shows the CCW duration distribution for $\alpha_1 = 0.1$ on the interval $[0, 60]$. 
        Here $\alpha_2 = -2$, $t_0 = 300$, $t_1 = 30$, $\bar{Y} = 5$, $T_m = 6000$ and $\sigma = 0.456$.} 
        \label{different_alpha_1}
\end{figure}

\begin{figure}[ht]
\centering
    \subfloat[]{\includegraphics
    [scale =0.5]{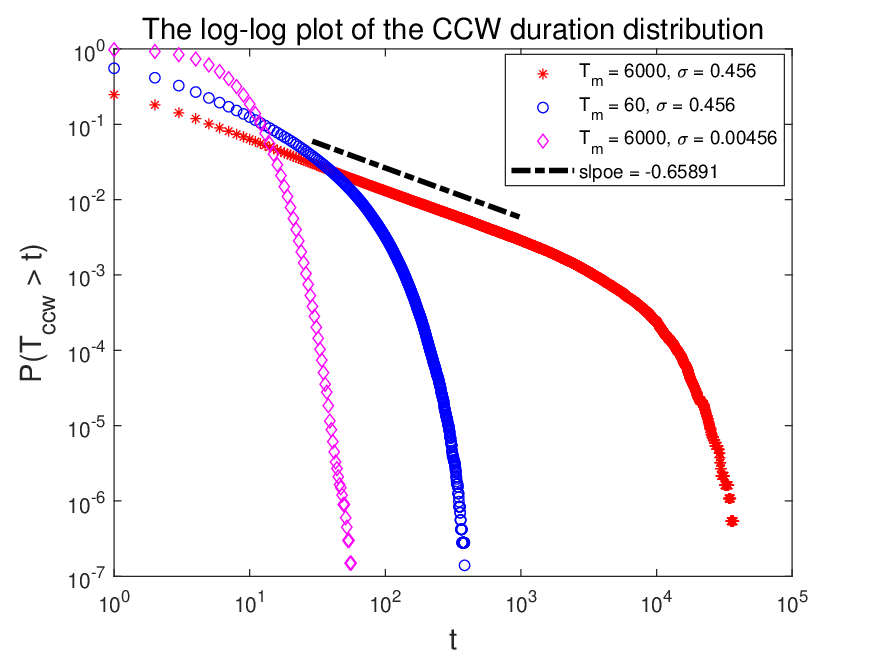}}
    \subfloat[]{\includegraphics
    [scale =0.5]{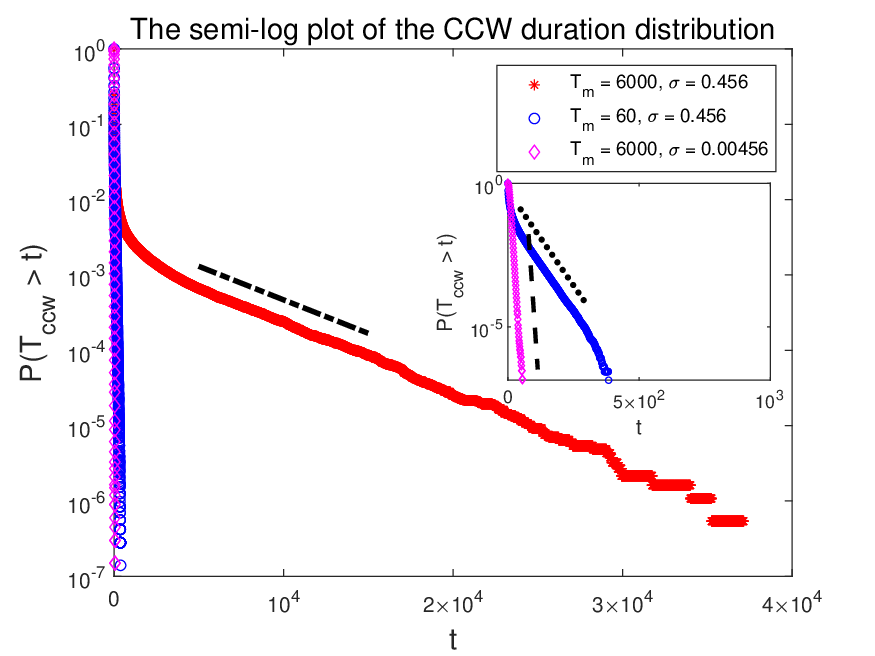}}
    \caption{The distribution of CCW duration for different values of $T_m$ and $\sigma$ in \eqref{mtequation}. (a) The log-log plot of $P(T_{ccw} > t)$ for $(T_m, \sigma)=(6000, 0.456)$ (red stars), $(60, 0.456)$ (blue circles) and $(6000, 0.00456)$ (magenta diamonds). 
    (b) The semi-log plot for $(T_m , \sigma)=(6000, 0.456)$ (red stars), $(60, 0.456)$ (blue circles) and $(6000, 0.00456)$ (magenta diamonds). The straight dash-dot line is fitted for the red stars between $5 \times 10^3$ and $1.5 \times 10^4$. The inset shows the CCW duration distribution for $(T_m = 60, \sigma = 0.456)$ and $(T_m = 6000, \sigma = 0.00456)$ on the interval $[0, 1000]$. 
    Here $\alpha_1 = 10$, $\alpha_2 = -2$, $t_0 = 300$, $t_1 = 30$ and $\bar{Y} = 5$.} 
    \label{Tdifferent_Tm_noise_sigma}
\end{figure}

An interesting question arises: 
what conditions must the parameters of the two-state model satisfy to achieve a power-law to exponential transition for the cumulative distribution function (CDF) of the CCW duration time? We have identified two necessary conditions: 
\begin{itemize}
\item[Cond I:] The motor response has to be ultrasensitive, i.e. $\alpha_1$ is large, so that $\Lambda_1(Y_t)$ transitions sharply near $\bar Y$. 
Two different values of $\alpha_1$ have been tested: the first one is the same as in \cite{tu2005white} with $\alpha_1 = 10$, while the second test is for $\alpha_1 = 0.1$ with all other parameters unchanged as described in Section \ref{two-statemodeltu}. As shown in Figure \ref{different_alpha_1}(a), when $\alpha_1=10$, the value of $\Lambda_1(Y_t) = \exp \left(\alpha_1 \frac{\Delta Y}{\bar Y}\right)$ increases rapidly from a small value to a large number near $\Delta Y = 0$, while the change is slower for $\alpha_1=0.1$. Figures \ref{different_alpha_1}(b) and \ref{different_alpha_1}(c) display the semi-log and log-log plots of the CDF of the CCW duration time, respectively. One can observe that when $\alpha_1 = 10$ the CDF exhibits desired transition in decay rate 
whereas when $\alpha_1=0.1$, the transition is not as evident. 

\item[Cond II:] The adaptation time $T_m$ and the noise intensity $\sigma$ must be large enough. As a consequence, the equilibrium distribution of $Y_t$ spreads out widely. 
Three sets of $(T_m,\sigma)$ are tested: $(6000,0.456)$, $(60,0.456)$, and $(6000,0.00456)$ with other parameters unchanged as described in Section \ref{two-statemodeltu}. The pair $(6000,0.456)$ is the same as in Section \ref{two-statemodeltu}. A shorter 
adaptation time is considered in $(60,0.456)$, while the noise intensity is smaller in $(6000,0.00456)$. The CDF of the CCW duration time for these three sets of parameters are displayed in Figure \ref{Tdifferent_Tm_noise_sigma}. As shown in Figures \ref{Tdifferent_Tm_noise_sigma}(a) and (b), the CDFs with $(60, 0.456)$ and $(6000,0.00456)$ exhibit exponential decay, 
while there is a transition from a power-law to exponential decay for the CDF when $T_m = 6000$ and $\sigma = 0.456$. 
Therefore, from Figure \ref{Tdifferent_Tm_noise_sigma}, the intermediate power-law decay can only be observed when both $T_m$ and $\sigma$ are large. 
\end{itemize}

\begin{figure}[!h]
\centering
    \subfloat[]{\includegraphics
    [scale =0.5]{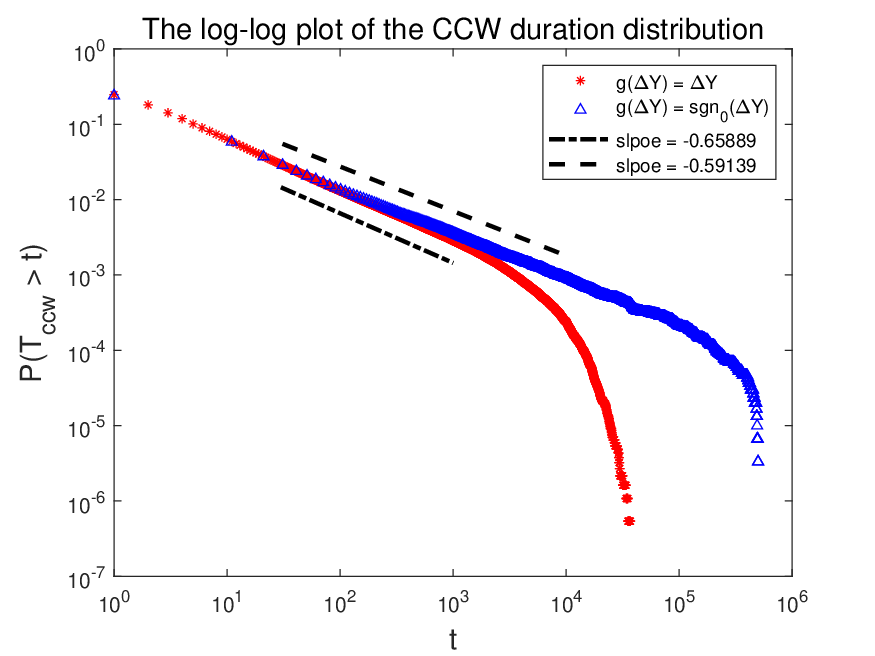}}
    \subfloat[]{\includegraphics
    [scale =0.5]{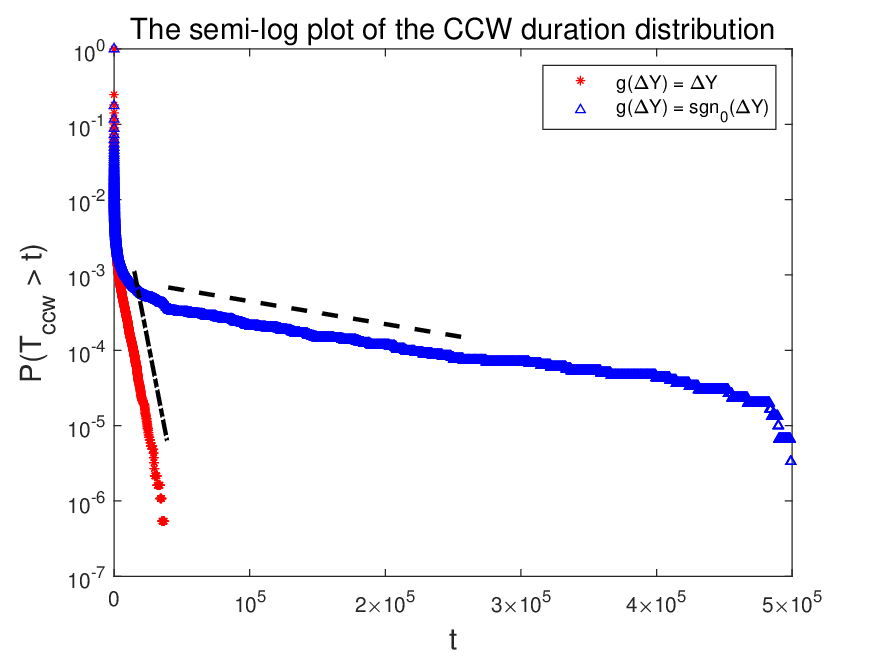}}
    \caption{The distribution of CCW duration for two different $g(\Delta Y)$. (a) The log-log plot of $P(T_{ccw} > t)$ for $g(\Delta Y) = \Delta Y$ (red stars) and $g(\Delta Y) = \text{sgn}_0(\Delta Y)$ 
    (blue triangles). 
    (b) The semi-log plot with a logarithmic y-axis for $P(T_{ccw} > t)$ with the two $g(\Delta Y)$. 
    Here, $\alpha_1 = 10$, $\alpha_2 = -2$, $t_0 = 300$, $t_1 = 30$ and $\bar{Y} = 5$ $T_m = 6000$ and $\sigma = 0.456$. }
    \label{Tdifferent_gm_sigma}
\end{figure}

\begin{figure}[!h]
\centering
    \includegraphics
    [scale =0.7]{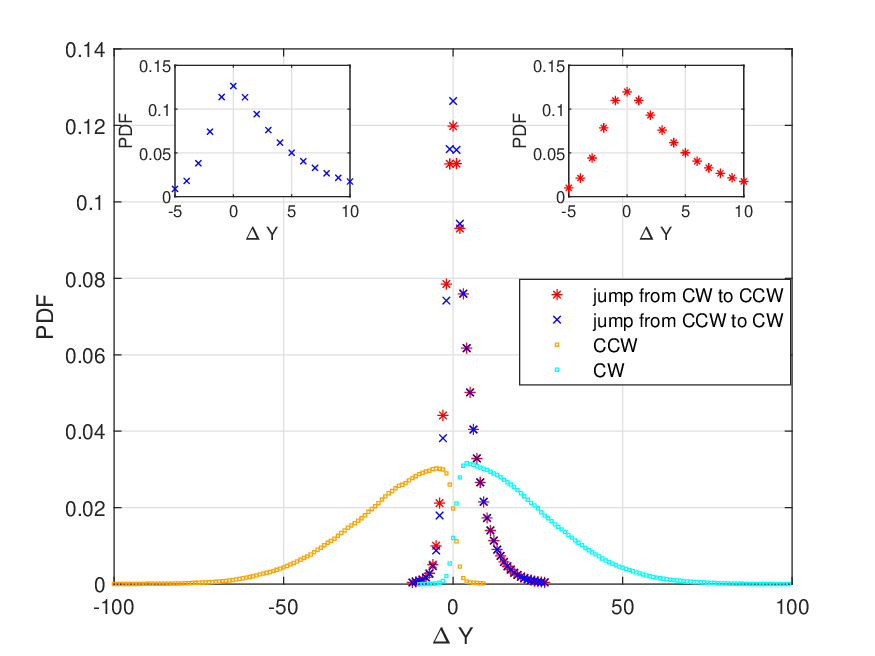}
    \caption{The distribution of $\Delta Y$ when the cells are at different rotational states. CCW state: orange squares; CW state: cyan squares. Particles from CW to CCW: red stars; from CCW to CW: blue crosses. The two insets display the distributions of red stars and blue crosses near the origin. Here  $\alpha_1 = 10$, $\alpha_2 = -2$, $t_0 = 300$, $t_1 = 30$, $\bar{Y} = 5$, $T_m = 6000$ and $\sigma = 0.456$. }
    \label{Ydifferent_distribution}
\end{figure}

In addition to Cond I and Cond II, 
we also make 
two observations: 
\begin{itemize}
\item[Ob I:]  Consider the following SDE which is a generalization of the OU process in \eqref{TupaperY}:
\begin{equation}\label{mtequation}
  d\Delta Y = -\frac{1}{T_m} g(\Delta Y) dt + \sigma dB_t.
\end{equation}
We observe that as long as $g(\Delta Y) \Delta Y >0$ for $\Delta Y \neq 0$, 
the specific form of $g(\Delta Y)$ does not seem essential. 
Indeed, two different forms of $g(\Delta Y)$ are tested. The first one is $g(\Delta Y) = \Delta Y$,
as used in 
\cites{tu2005white, XZZT}. The second one is a step function
\begin{equation}\label{gm}
    g(\Delta Y) = \text{sgn}_0(\Delta Y) := \begin{cases}
                1, \quad \Delta Y \geq 0, \\
                -1, \quad \Delta Y < 0.
           \end{cases}
\end{equation}
In both cases, the first term on the right-hand side of \eqref{mtequation} drives $\Delta Y$ 
toward zero. 
In Figure \ref{Tdifferent_gm_sigma}, we display the CCW duration distributions for the two choices of $g(\Delta Y)$ with the same parameter as in Section \ref{two-statemodeltu}. 
It is clear that the CDFs of both cases exhibit a crossover from an intermediate power-law 
to a subsequent exponential decay. 
\item[Ob II:] In Figure \ref{Ydifferent_distribution}, we use the same parameters as in Section \ref{two-statemodeltu}. Figure \ref{Ydifferent_distribution} shows the distribution of $\Delta Y$ when the cells are at different rotating states. It can be seen that 
most CCW (CW) rotation occurs when $\Delta Y < 0$ ($\Delta Y > 0$). 
More specifically, when cells are at the CCW state, $\Delta Y$ is mainly distributed over the interval $[-100, 10]$, and when at the CW state, $\Delta Y$ is distributed over the interval $[-10, 100]$. 
Moreover, most jumps from CCW to CW (or CW to CCW) happen around $\Delta Y = 0$.
\end{itemize}

Cond I, Cond II, Ob I, and Ob II lead us to the construction of a simplified one-state model in the next subsection, including the choices of the scaling regime, the tumbling rate $\Lambda$ and the adaptation function $g$.

\subsubsection{The simplified one-state model}\label{twostatemodel}

\begin{figure}[!h]
\centering
    \subfloat[]{\includegraphics
    [scale =0.5]{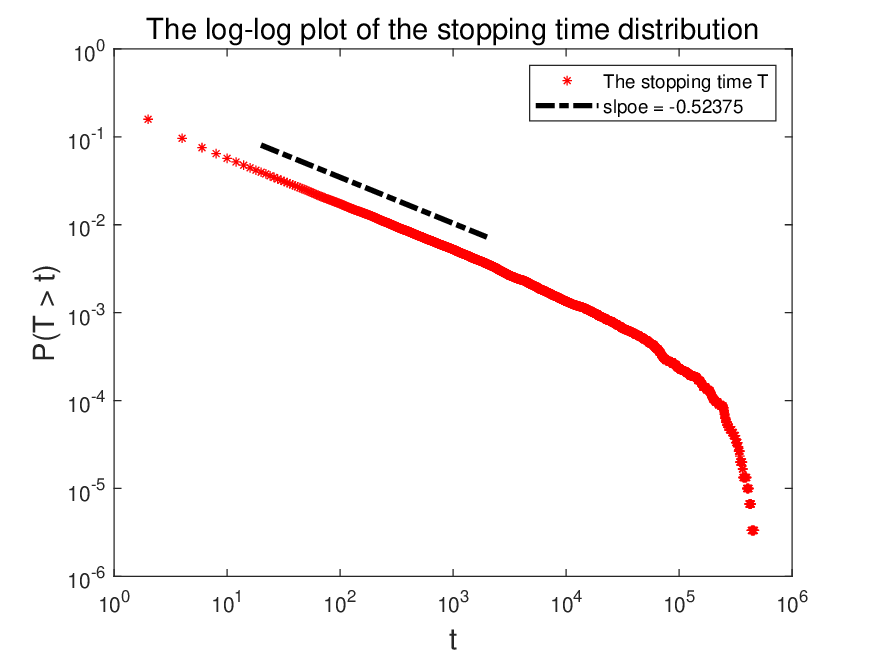}}
    \subfloat[]{\includegraphics
    [scale =0.5]{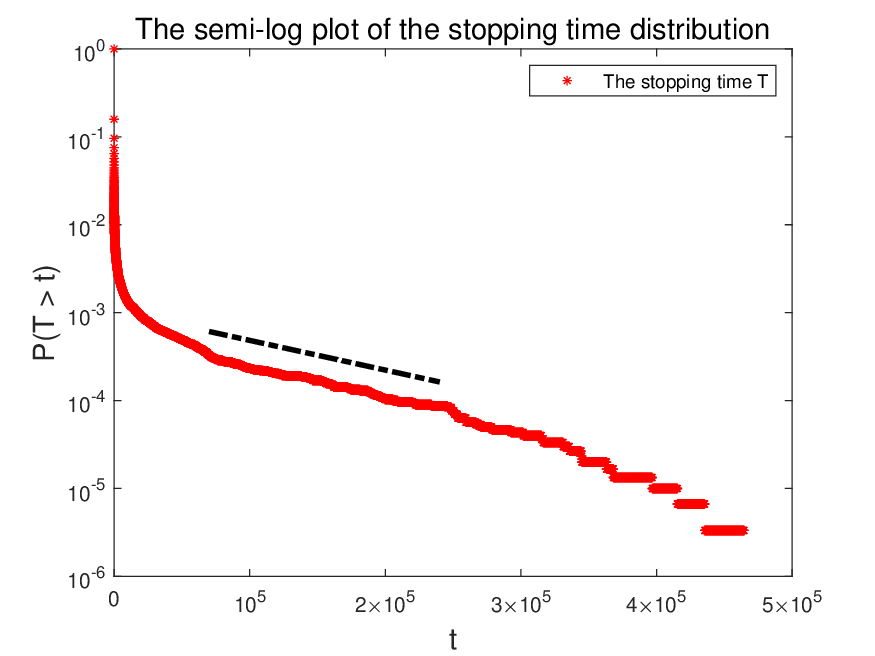}}
    \caption{The distribution of the stopping time $T$ for the one-state model with $T_m = 6000$ and $\sigma = 0.456$. (a) The log-log plot for $P(T > t)$. The slope $-0.57362$ is fitted over the interval $[10^2, 10^4]$. (b) The semi-log plot with a logarithmic y-axis for $P(T > t)$. The straight dash-dot line is fitted over the interval $[6 \times 10^4, 2.6 \times 10^5]$.}
    \label{one-state model fig}
\end{figure}

Many models in the literature for \textit{E. coli} chemotaxis ignore the tumbling time, i.e., the CW rotational state \cites{taktikos2013motility, Xue2021, levywalks}. These simplified models are easier to analyze and sufficient to explain many experimental observations. For simplicity, in this paper we only consider the CCW state, which results in a one-state model.

Let $m_t = \Delta Y$. We consider the scenario that the particle jumps outside the CCW state and immediately starts another CCW state according to a state-dependent Poisson clock with a rate $\Lambda(m_t)$. We choose
\begin{equation}\label{Lambda}
    \Lambda(m_t) = \begin{cases}
                 1, & \text{$m_t \geq 0$}, \\
                 0, & \text{$m_t < 0$}.
            \end{cases}
\end{equation}
Thus, $\Lambda(m_t)$ increases immediately from zero to an $O(1)$ constant at $m_t = 0$, which mimics Cond I, i.e., the ultra-sensitivity of the motor response.
The associated stopping time can be calculated as 
\begin{equation}\label{stoppingtime}
    T = \inf \left\{t: \int_{0}^{t} \Lambda(m_s) \, ds > \Gamma \right\},
\end{equation}
where $\Gamma \sim \exp(1)$ is an exponentially distributed random variable with rate $1$ and is independent of $m_t$. To verify that the stopping time has the form as in \eqref{stoppingtime}, one can denote $\mathcal{F}_t$ as the $\sigma$-algebra generated by $m_s$ for $s \le t$. Then
by \cite{XZZT}, we have
\begin{equation*}
    \mathbb{P}(T > t | \mathcal{F}_t) = \mathbb{P}\left( \left. \int_{0}^{t} \Lambda(m_s) \, ds \leq 
    \Gamma \right| \mathcal{F}_t\right) = \exp\left(-\int_{0}^{t} \Lambda(m_s) \, ds \right).
\end{equation*}
Given $\mathcal{F}_t$, the conditional jumping rate at time $t$ is then given by
\begin{equation*}
     \frac{\frac{d \mathbb{P}(T \leq t | \mathcal{F}_t)}{dt}}{\mathbb{P}(T > t | \mathcal{F}_t)} = \frac{\Lambda(m_t) \exp\left(-\int_0^t \Lambda(m_s) \, ds\right)}{\exp\left(-\int_0^t \Lambda(m_s) \, ds\right)} = \Lambda(m_t).
\end{equation*}

According to Ob II, one can assume that when particles jump outside of the CCW state, $m_t$ is immediately reset to the value $0$. This results in the following one-state model:
\begin{subnumcases}{\label{one-statemodel}}
    dm_t = -\frac{1}{T_m} g(m_t) \, dt + \sigma dB_t, \\
    T = \inf \left\{t: \int_{0}^{t} \Lambda(m_s) \, ds > \Gamma \right\}, \\
    m_t |_{t=T+} = 0,
\end{subnumcases}
where $g(m_t)$ and $\Lambda(m_t)$ are defined in \eqref{gm} and \eqref{Lambda} respectively.
We simulate the one-state model defined in \eqref{one-statemodel} with $T_m = 6000$ and $\sigma = 0.456$ as in Section \ref{two-statemodeltu}.
There are $N = 100$ samples tracked, and the total running time is $6 \times 10^5$. Each sample is represented by its internal state $m_t$ and the stopping time $T$. Details of the numerical algorithm are described in Appendix \ref{one-statemodelappendix}. Figure \ref{one-state model fig} shows the CDF of the stopping time $T$ of the one-state model \eqref{one-statemodel}. One can observe that the CDF of the stopping time $T$
exhibits an intermediate power-law decay and a subsequent transition to exponential decay.

\subsection{Kinetic model with cell movement} \label{sec:kinetic}
In this subsection, we introduce an individual-based model that integrates the run-and-tumble motion with the internal adaptation dynamics modeled by~\eqref{one-statemodel}. 

\subsubsection{Kinetic model with internal dynamics}\label{IBM}
To consider the dispersion of particles, we incorporate the particle movements into the one-state model~\eqref{one-statemodel}. Each particle is characterized by its position $x \in \mathbb{R}^d$, velocity direction $v \in \mathbb{S}^{d-1}$ (a unit spherical surface in $\mathbb{R}^d$), and internal state $m \in \mathbb{R}$. Between two successive jumps, the particle evolves as
\begin{subnumcases}{\label{IBMmodel}}
    dx = V_0v dt, \label{ibmnew1} \\
    d v = 0, \label{ibmnew2} \\
    dm_t = -\frac{1}{T_m} g(m_t) dt + \sigma dB_t, \label{ibmnew3}
\end{subnumcases}
where $V_0$ is a constant representing the magnitude of the cell movement velocity. The stopping time $T$ is defined in \eqref{stoppingtime}. At time $T$, the internal state $m_t$ is reset to $0$, and the particle randomly chooses a new movement direction with uniform probability. More precisely,
\begin{equation}\label{IBMmodel2}
    \begin{cases}
        T = \inf \left\{t:\int_{0}^{t} \Lambda(m_s) ds > \Gamma \right\}, \\
        m|_{t = T+} = 0, \\
        x|_{t = T_+} = x|_{t = T_-}, \\
        v|_{t = T_+} \sim U(\mathbb{S}^{d-1}),
        \end{cases}
\end{equation}
where $\Gamma \sim \exp(1)$ and $U(\mathbb{S}^{d-1})$ denotes the uniform distribution on the sphere $\mathbb{S}^{d-1}$.

\subsubsection{The crossover from ballistic transport to normal diffusion}\label{sec:222}
We conduct numerical simulations of the individual-based model proposed in \eqref{IBMmodel} and \eqref{IBMmodel2}. In the simulation, let the superscript $n$ be the index for the cell. Each cell evolves as a particle with position $x^n$, velocity direction $v^n$, and internal state $m^n$. We evolve $(x^n, v^n, m^n)$ according to 
\eqref{IBMmodel} and \eqref{IBMmodel2}. One spatial dimension is considered, and the movement velocity direction $v^n \in \{-1, 1\}$. The initial values $x^n$, $m^n$ for all samples are $0$, and the initial $v^n$ is randomly chosen to be $1$ or $-1$ with equal probability. We run $N = 10^3$ samples, and the total running time $T_t = 10^6$ for each sampled particle. Let $V_0 = 0.02$. 
We use $\sigma = O(1)$ and check the effects of different values of $T_m$ on the macroscopic dispersion behavior of the particles. Here, to reduce the complexity of determining the asymptotic regimes, we fix $\sigma = \sqrt{2}$ 
and choose $T_m = 10^{-2}, 1, 10, 100, 1000$, and $+\infty $. When $T_m = \infty$, or equivalently, $\frac{1}{T_m} = 0$, the drift term in \eqref{ibmnew3} is neglected. The time step is chosen to be $\Delta \tau = 1$ for $T_m = 1, 10, 100, 1000$, and $\infty$. When $T_m = 10^{-2}$, we have $\frac{1}{T_m} = 100$, which implies a large drift term in \eqref{ibmnew3}. To reduce errors and ensure the accuracy of numerical simulations, we use a smaller $\Delta \tau = 10^{-2}$ and the run time is also $T_t = 10^6$. The details of the numerical algorithms are in Appendix \ref{IBMappendix}. 

In the subsequent part, we use three different criteria to verify the crossover of different dispersion patterns: the mean square displacement (MSD), the run time distribution (RTD), and the PDF of $\Delta x/\Delta t^\beta$. Here, $\Delta t$ is a fixed time interval, and $\Delta x$ is the particle displacement during this time interval. Furthermore, according to \cites{levywalks, yuanliPNAS, SwimmingEcoli, PhysRevE.91.022131}, if in a given time interval ${\rm MSD} \sim t^{\beta_0}$ with $\beta_0\in(1,2)$ being a constant, the PDF of the rescaled displacement $\Delta x/\Delta t^{\beta}$ with $\beta = 1/(3-\beta_0)$ will overlap for different values of lag time $\Delta t$ \cites{yuanliPNAS}. The rescaled displacement is more localized in time and thus can be used to verify the crossover of different dispersion patterns with time.
Among these three measurements, the MSD is on the population level, while the other two are more for the PDFs of individuals\cite{yuanliPNAS}.

\begin{figure}[ht]
\centering
\subfloat[]{\includegraphics
    [scale =0.5]{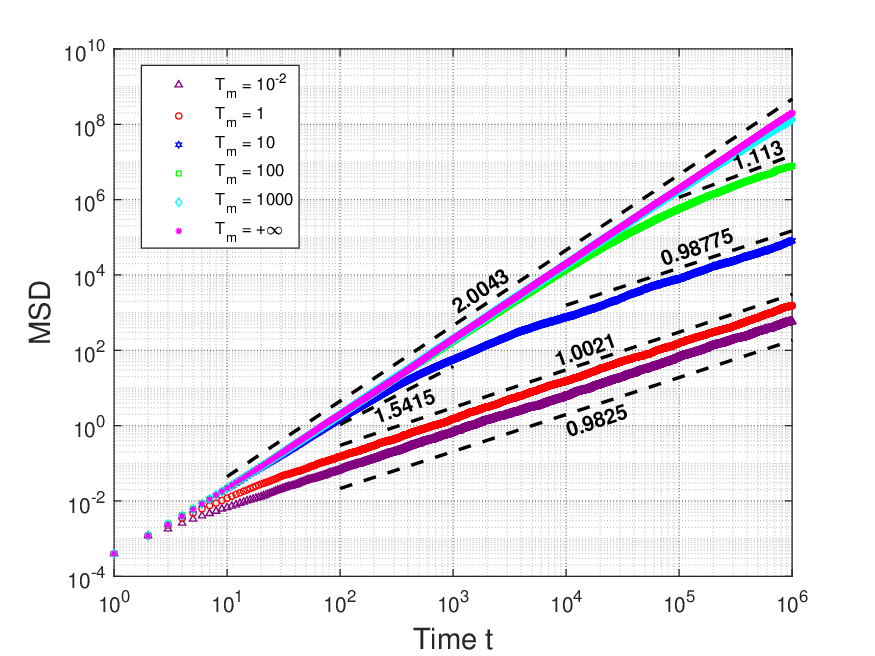}}
    \subfloat[]{\includegraphics
     [scale =0.5]{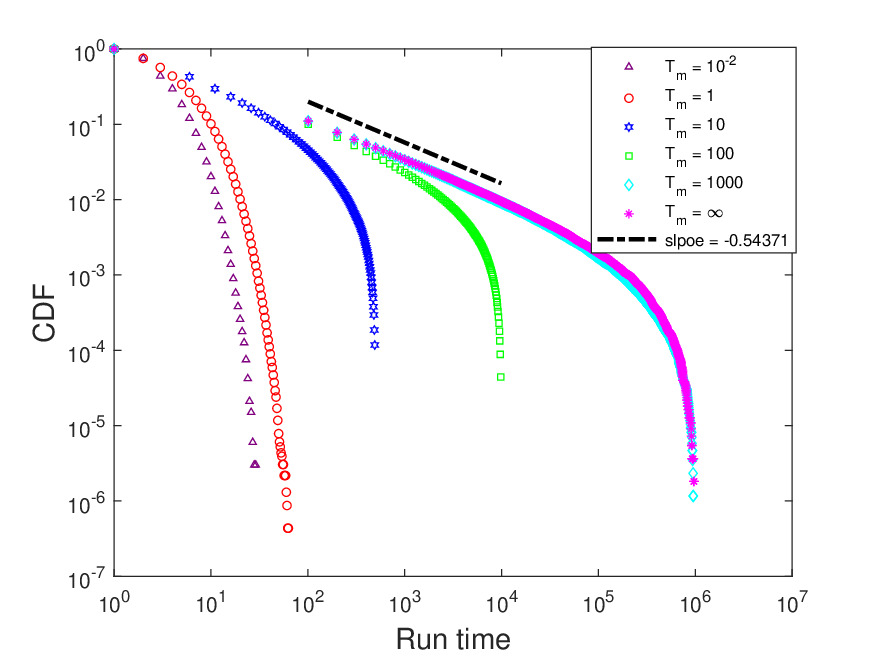}}
    \caption{(a) The log-log plot of the MSD for different $T_m$. $ T_m = 10^{-2}$ : purple triangles; $ T_m = 1$ : red circles; $ T_m = 10$ : blue hexagrams; $ T_m = 100$ : green squares; $ T_m = 1000$ : cyan diamonds; $ T_m = +\infty$ : magenta stars. (b) The log-log plot of the RTD for different $T_m$. $ T_m = 10^{-2}$ : purple triangles; $ T_m = 1$ : red circles; $ T_m = 10$ : blue hexagrams; $ T_m = 100$ : green squares; $ T_m = 1000$ : cyan diamonds; $ T_m = +\infty$ : magenta stars.} 
    \label{MSD}
\end{figure}

\paragraph{MSD.} The MSD is defined by
\begin{equation*}
    \mathrm{MSD} := \left \langle  (x-x_0)^2 \right \rangle = \frac{1}{N} \sum_{n=1}^N (x^n(t)-x^n(0))^2,
\end{equation*}
where $x^n(t)$ is the position of the $ n $-th particle at time $ t $, $ x^n(0) $ is the initial position of the $ n $-th particle, and $ \langle \cdot \rangle $ denotes the average over all sampled particles.
We plot the MSD for different $ T_m = 10^{-2}, 1, 10, 100, 1000$, and $+\infty$ in Figure \ref{MSD}(a). For $ T_m = 10^{-2}$ and $1$, $ \mathrm{MSD} \sim t $ in the time interval $[10^2, 10^6]$. This indicates that particles undergo normal diffusion at the macroscopic scale. When $ T_m = 10 $ and $ 100 $, $ \mathrm{MSD} \sim t^2 $ in the interval $[10 , 10^2 ]$ and $ \mathrm{MSD} \sim t $ when $ t \in [10^5 , 10^6 ]$, which implies a transition from ballistic transport to normal diffusion at the population level. As for $ T_m = 1000  $ and $ T_m = +\infty $, it can be observed that $ \mathrm{MSD} \sim t^2 $ in the interval $[10 , 10^6]$, which indicates that particles manifest ballistic transport. To summarize, when the adaptation time $ T_m $ is short ($ T_m = 10^{-2} , 1 $), particles display normal diffusion, while when $ T_m $ is long ($ T_m = 10, 100 $), particles initially exhibit ballistic transport and then shift to normal diffusion. When $ T_m $ is large enough ($ T_m = 1000, +\infty $), particle motion presents a ballistic behavior throughout the entire simulation time $ T_t = 10^6 $. Furthermore, as the adaptation time increases, the crossover time at which ballistic transport transitions to normal diffusion also increases. 

\begin{figure}[ht]
\centering
    \subfloat[]{\includegraphics
    [scale =0.36]{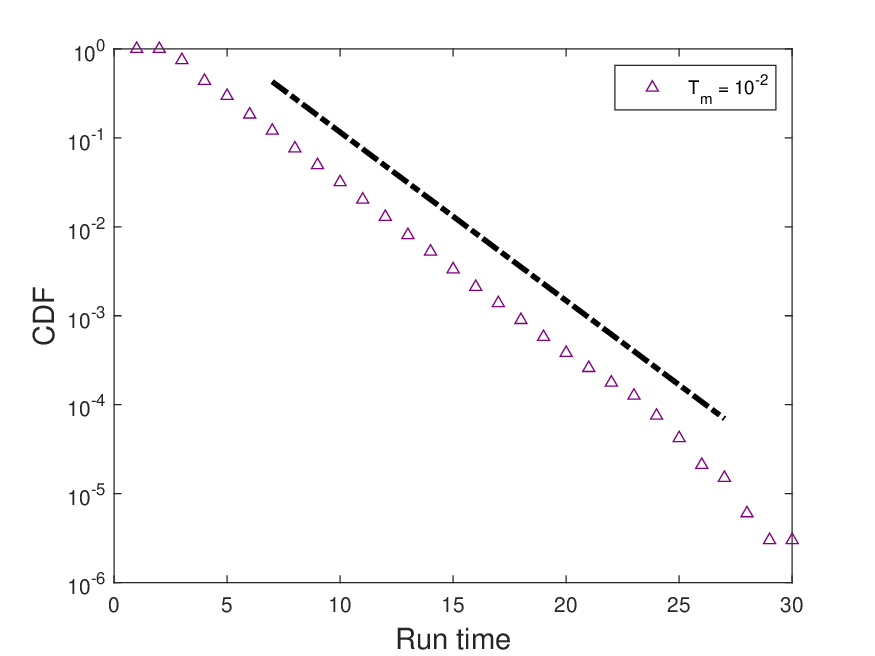}}
    \subfloat[]{\includegraphics
    [scale =0.36]{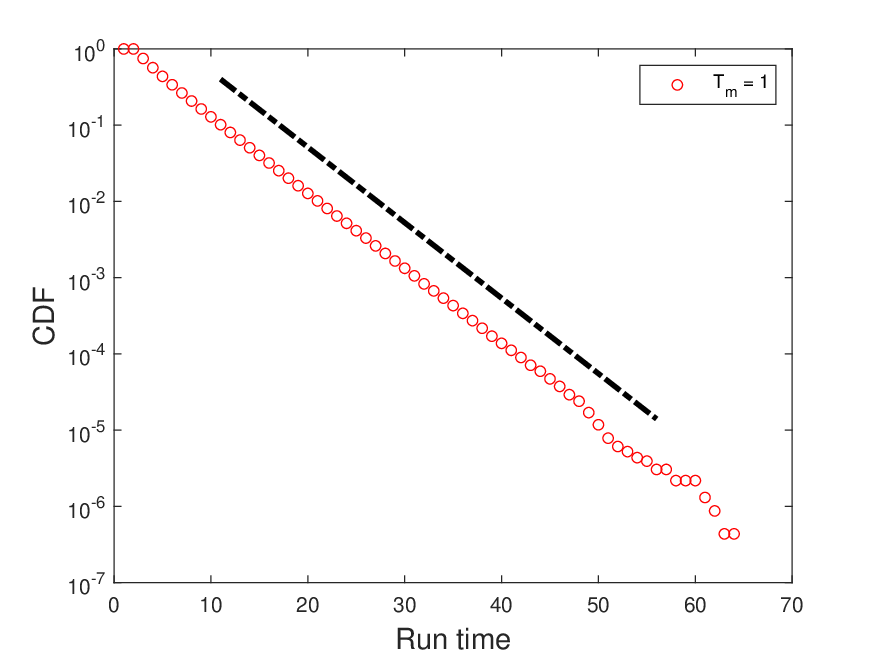}}
    \subfloat[]{\includegraphics
    [scale =0.36]{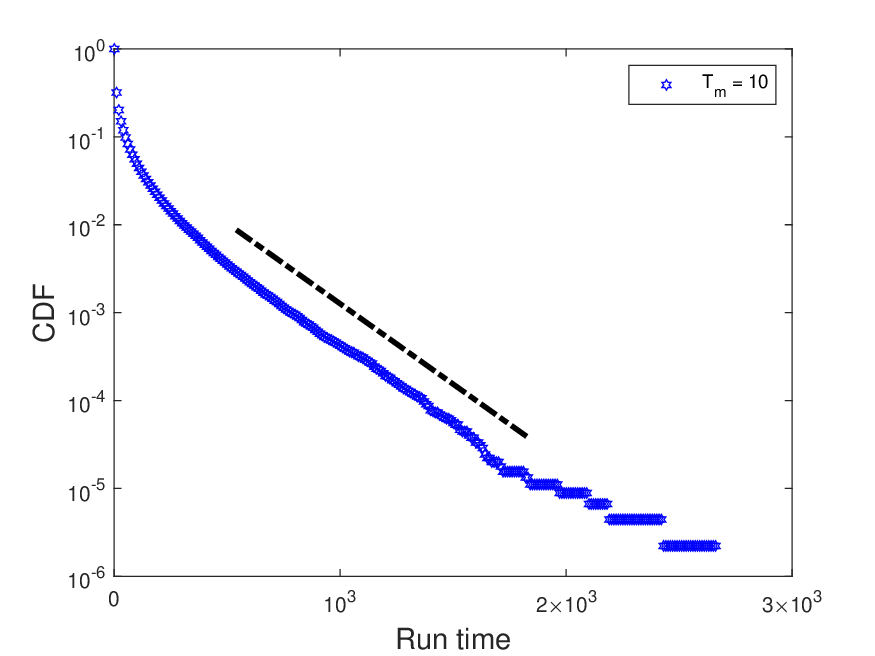}}
    \\
    \subfloat[]{\includegraphics
    [scale =0.36]{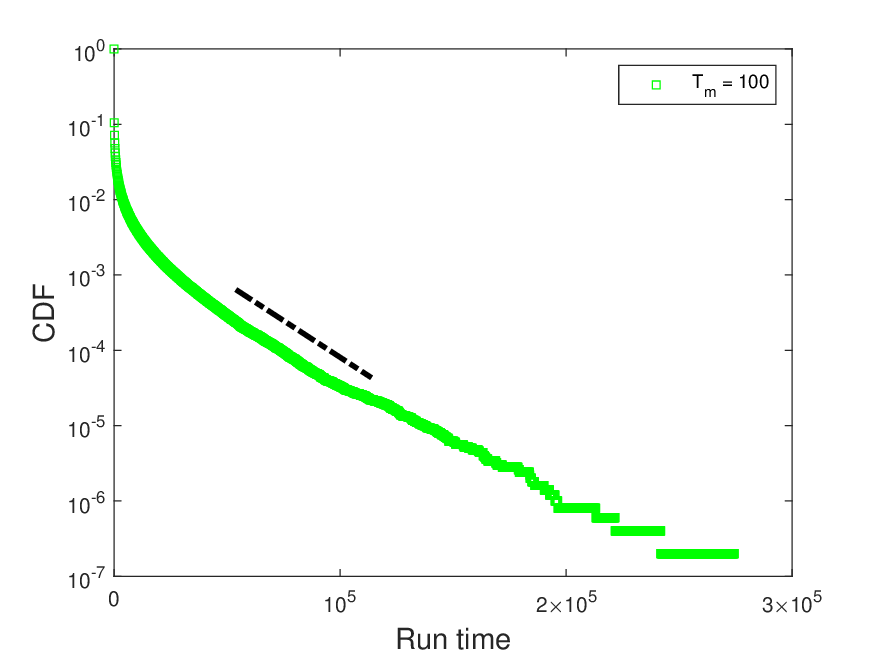}}
    \subfloat[]{\includegraphics
    [scale =0.36]{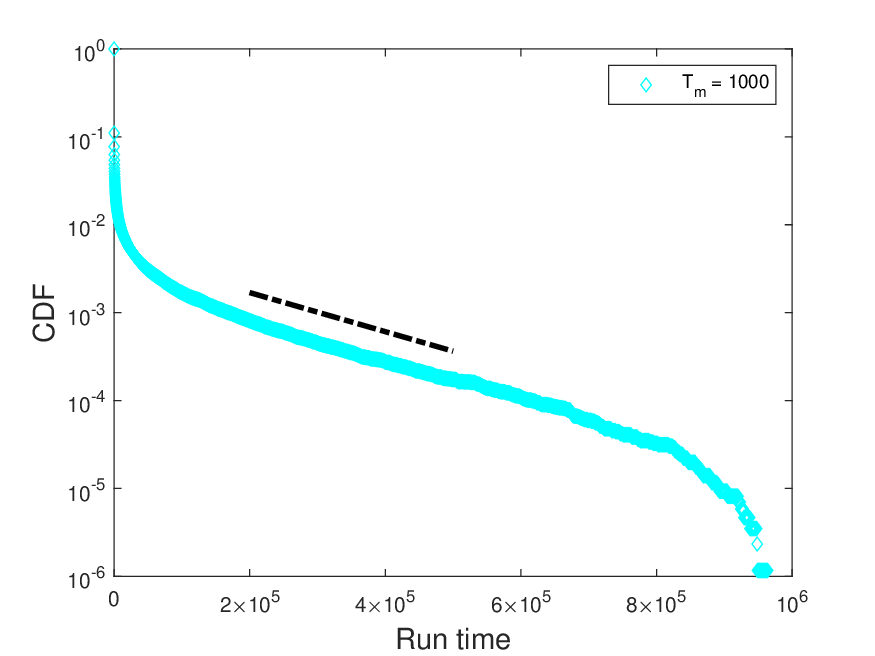}}
    \subfloat[]{\includegraphics
    [scale =0.36]{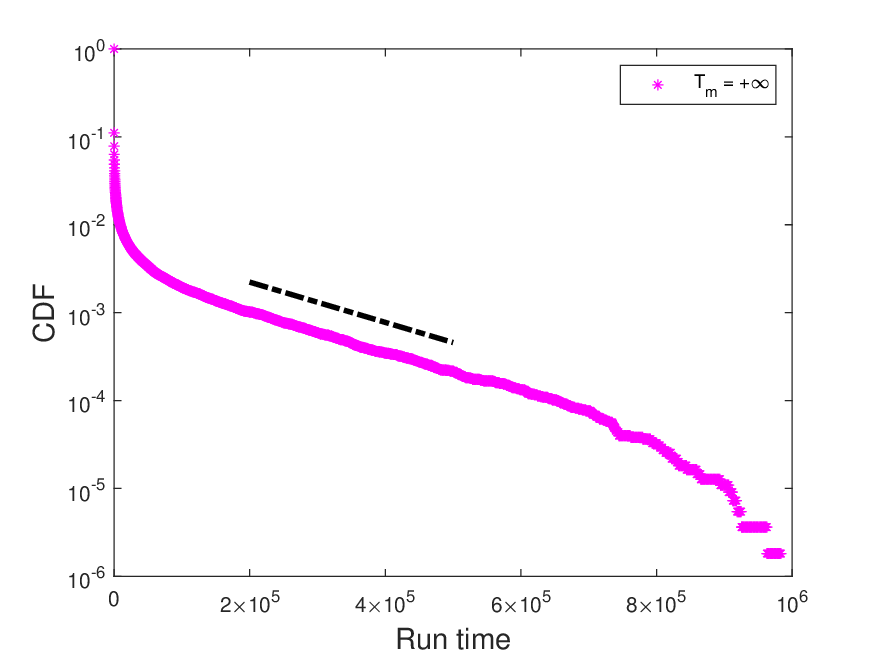}}
    \caption{The semi-log plots of RTD for different $T_m$. (a) $T_m = 10^{-2}$; 
    (b) $T_m = 1$; (c) $T_m = 10$; (d) $T_m = 100$; (e) $T_m = 1000$; (f) $T_m = +\infty$.}
    \label{RLDexp}
\end{figure}

\paragraph{The run time distribution (RTD).}
We stop the simulation at time $T_t$ and collect all run times $\{T_{\text{runtime}}\}$, which are measured between two consecutive changes in movement directions for each of the $N = 10^3$ particles.
Since the magnitude of the velocity $V_0$ is constant, the run time is proportional to the run length. Figure \ref{MSD}(b) shows the log-log plot of the CDF of the RTD, defined as $P(T_{\text{runtime}} > \text{Run time})$. From Figure \ref{MSD}(b), for $T_m = 100, 1000$, and $+\infty$, there exists an interval where the RTD exhibits a power-law decay. Moreover, all curves in Figure \ref{MSD}(b) have a rapidly decaying tail.

To delve deeper into the rapidly decaying tail in Figure \ref{MSD}(b), semi-log plots of the RTD for all different $T_m$ are shown in Figure \ref{RLDexp}. We observe that for all $T_m$, the tail parts can be fitted by straight lines, indicating that the RTD for all $T_m$ decays exponentially fast at the tail end. The exponential decay tails may be due to the limited number of samples or the effect of the adaptation term, as discussed in \cite{XZZT}. 


\begin{figure}[ht]
\centering
\subfloat[]{\includegraphics
    [scale =0.36]{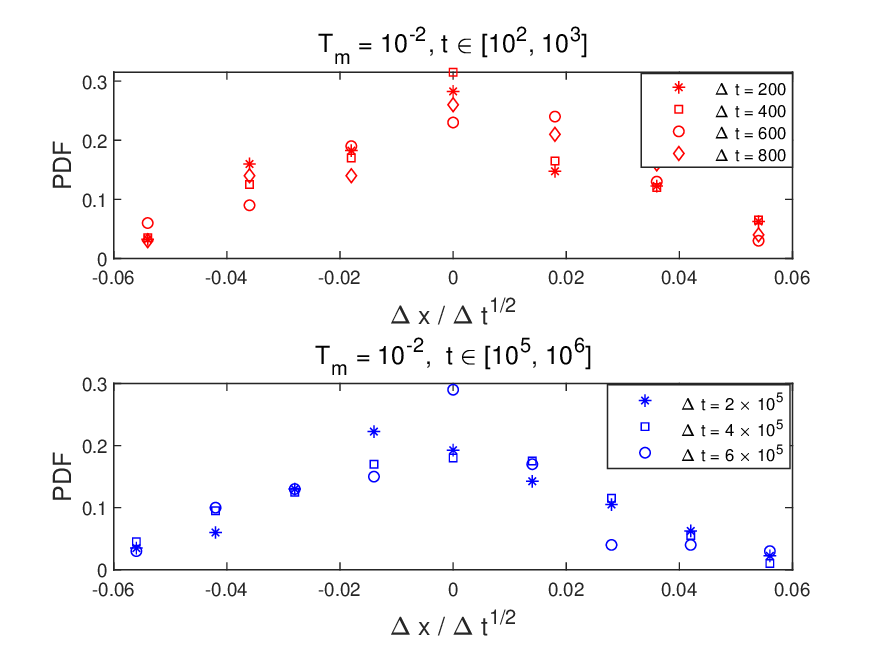}}
    \subfloat[]{\includegraphics
    [scale =0.36]{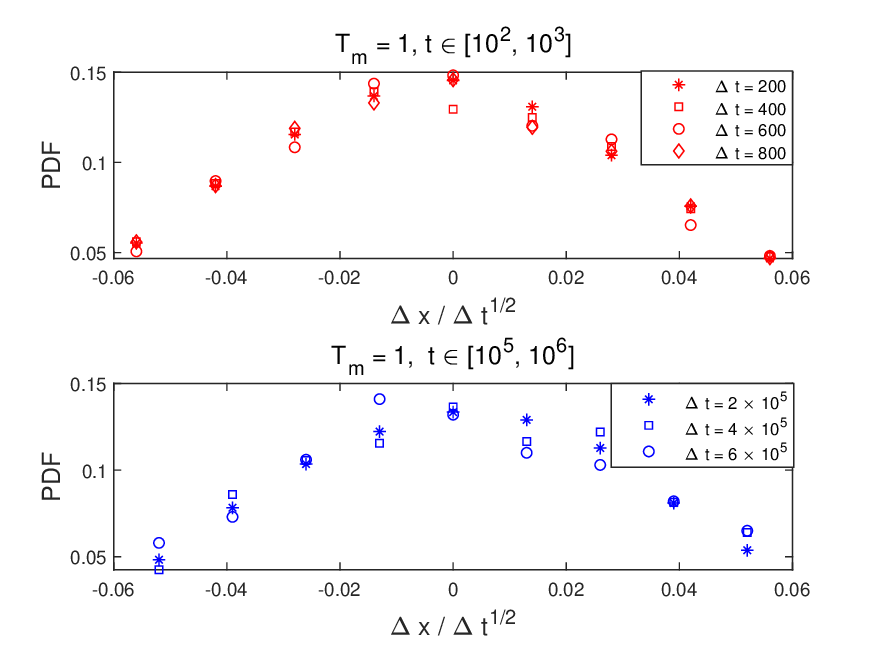}}
    \subfloat[]{\includegraphics
    [scale =0.36]{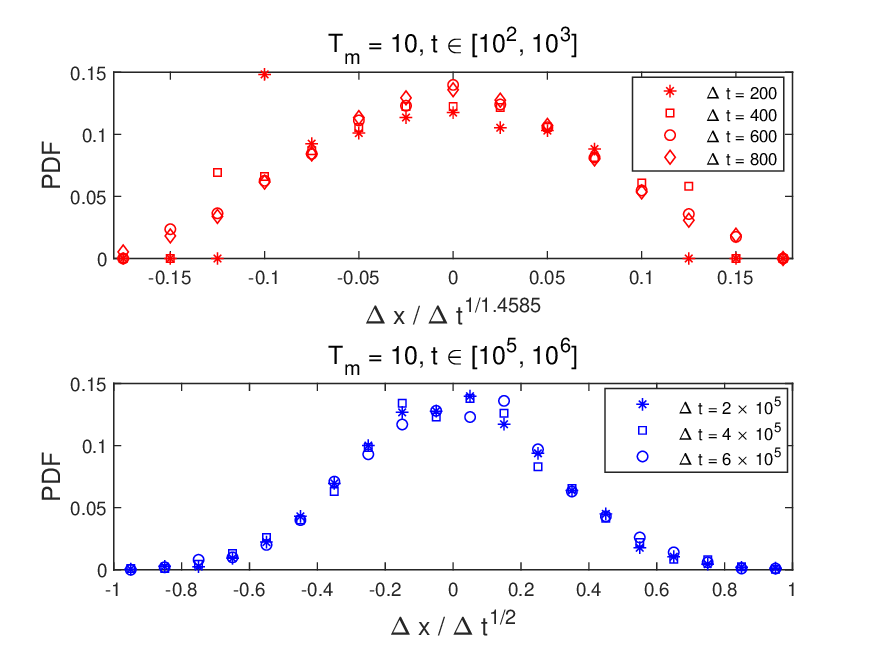}}
    \\
    \subfloat[]{\includegraphics
    [scale =0.36]{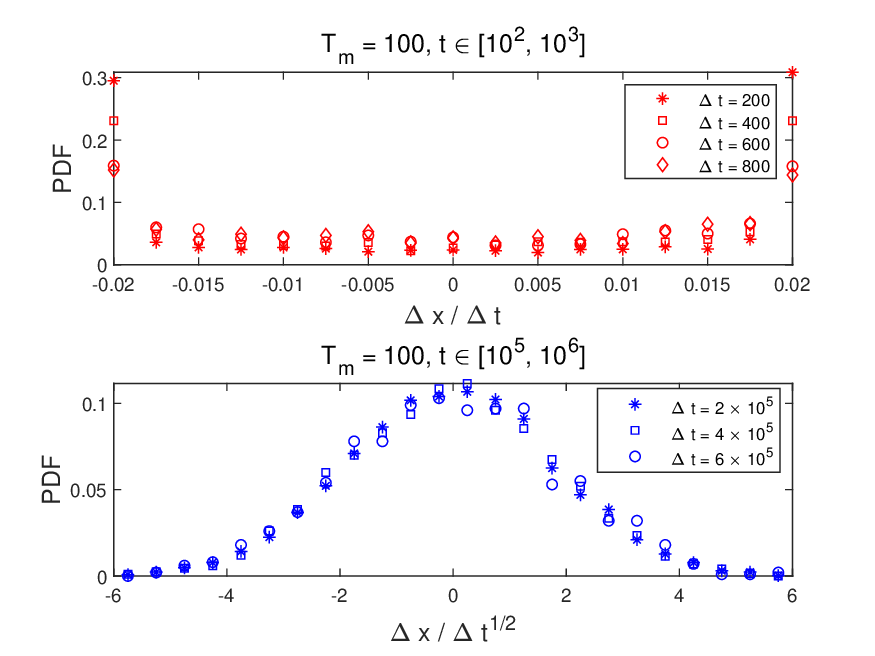}}
    \subfloat[]{\includegraphics
    [scale =0.36]{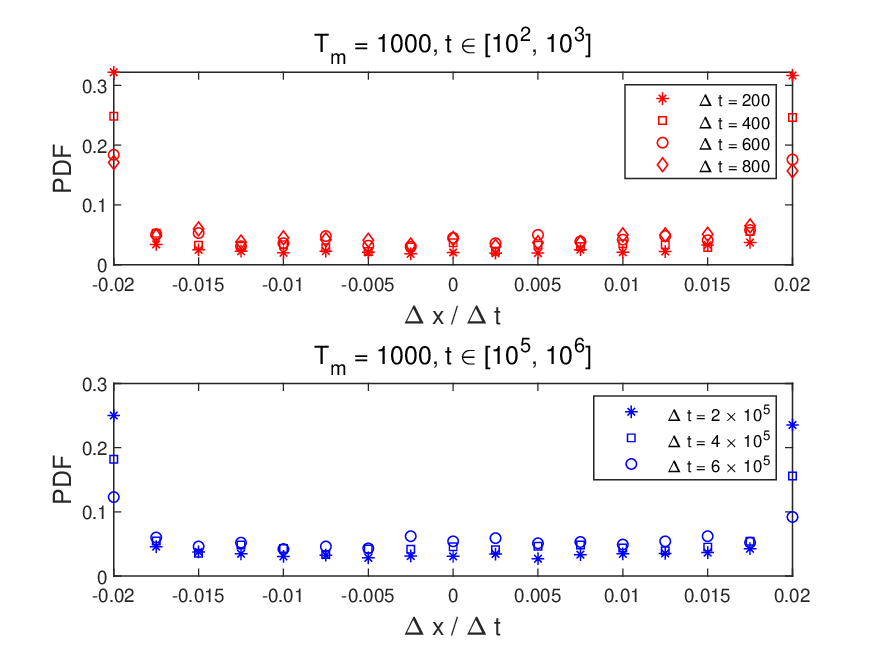}}
    \subfloat[]{\includegraphics
    [scale =0.36]{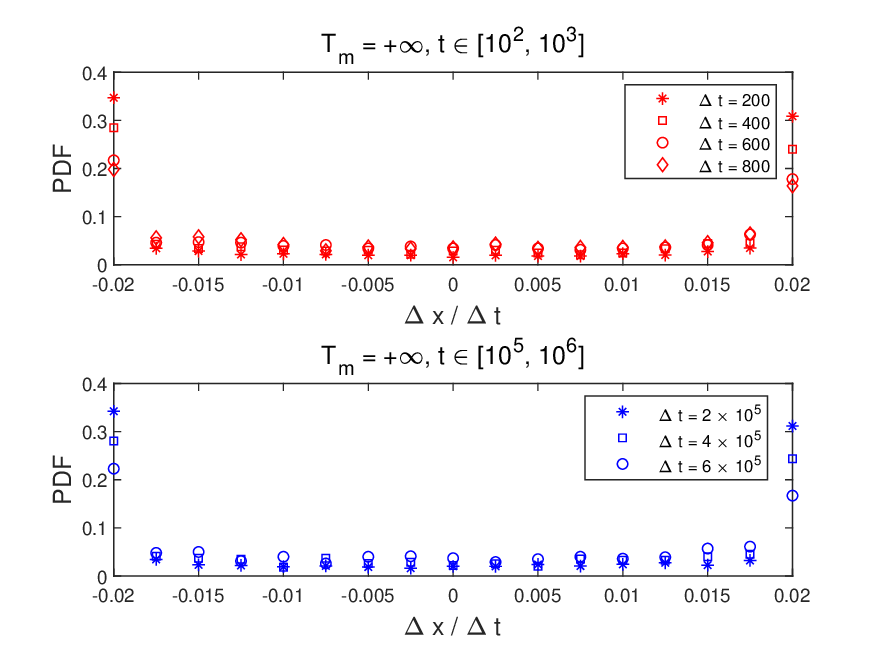}}
    \caption{The PDF of $\frac{\Delta x}{\Delta t^{\beta}}$ for different $T_m$ values and different $\Delta t$. (a) $T_m = 10^{-2}$, (b) $T_m = 1$, (c) $T_m = 10$, (d) $T_m = 100$, (e) $T_m = 1000$, (f) $T_m = +\infty$. In each of the subplots, the picture on the top is for $t \in [10^2, 10^3]$ where the PDFs for $\Delta t = 200, 400, 600, 800$ are respectively plotted with stars, squares, circles, and diamonds. The picture on the bottom of each subplot is for $t \in [10^5, 10^6]$, where the PDFs for $\Delta t = 2 \times 10^5, 4 \times 10^5, 6 \times 10^5, 8 \times 10^5$ are respectively plotted with stars, squares, and circles.}
    \label{dx_dt}
\end{figure}

\paragraph{PDF of $\Delta x / \Delta t^{\beta}$.}
To further validate the particle's dispersion pattern, we present the PDF of $\Delta x / \Delta t^{\beta}$ in Figure \ref{dx_dt}. For different dispersion patterns, the rescaled displacements take different values of $\beta$. For ballistic transport, $\beta = 1$, and for normal diffusion, $\beta = \frac{1}{2}$ as in \cites{levywalks, yuanliPNAS}. 
In each subplot of Figures \ref{dx_dt}, the top panel is for the time interval $[10^2 , 10^3]$ and the PDFs of 
\begin{equation*}
    \Delta x / \Delta t^{\beta} =   \left[x(10^2 + \Delta t) - x(10^2) \right]/ \Delta t^{\beta},
\end{equation*}
for different $\Delta t = 200 , 400 , 600 , 800 $ are displayed. From Figure \ref{MSD}(a), when $T_m = 10^{-2}$ and $1$, MSD $\sim t$. Therefore $\beta$ is set to $\frac{1}{2}$. When $T_m = 10$, we have MSD $\sim t^{1.5415}$. In this case, $\beta = 1/(3-1.5415) = 1/1.4585$. For $T_m = 100, 1000$ and $+\infty$, since MSD $\sim t^2$, we set $\beta = 1$.
%
The bottom panel of each subplot is for the time interval $[10^5, 10^6]$. The PDFs of 
\begin{equation*}
    \Delta x / \Delta t^{\beta} =  \left[x(10^5 + \Delta t) - x(10^5)\right] / \Delta t^{\beta},
\end{equation*}
are displayed for different $\Delta t = 2 \times 10^5 , 4 \times 10^5 , 6 \times 10^5$. From Figure \ref{MSD}(a), for $T_m = 10^{-2}, 1, 10$ and $100$, we always have MSD $\sim t$. Thus $\beta = \frac{1}{2}$ is chosen. For $T_m = 1000$ and $+\infty$, since MSD $\sim t^2$, we set $\beta = 1$. 
%
When $T_m$ is short ($T_m = 10^{-2}, 1$), particles exhibit normal diffusion, while when $T_m$ is long ($T_m = 100$), particles initially exhibit ballistic transport and transition to normal diffusion. When $T_m$ is sufficiently long ($T_m = 1000$ and $+\infty$), particles exhibit only ballistic transport. 
Additionally, for $T_m = 10$, Figure \ref{dx_dt}(c) shows the PDFs of $\Delta x / \Delta t^{\frac{1}{1.4585}}$ in $[10^2, 10^3]$ and $\Delta x / \Delta t^{\frac{1}{2}}$ when $t \in [10^5, 10^6]$. It can be observed that the PDFs for different $\Delta t$ do not coincide with each other, thus it is hard to identify the particular dispersion pattern in the interval $[10^2, 10^3]$. On the other hand, normal diffusion can be observed in the interval $[10^5, 10^6]$.


\section{The PDE model and main results} \label{sec:results}
\subsection{The PDE model and nondimensionalization}
The PDE model associated with the individual-based model \eqref{IBMmodel}--\eqref{IBMmodel2} describes the time evolution of the PDF $f(t,x,v,m)$ of the particles at time $t$, position $x \in \mathbb{R}^d$, velocity $v \in \mathbb{V}_0$, and internal state $m \in \mathbb{R}$. The velocity space $\mathbb{V}_0$ is $\{-V_0,V_0\}$ for $d=1$ and $\partial B(0,V_0)$ for $d \geq 2$ (The spherical surface with radius $V_0$). We denote $||\mathbb{V}_0|| = \int_{\mathbb{V}_0} \, dv$ as the surface area of the sphere $\partial B(0, V_0)$. 

By Ito's formula, the Fokker-Planck operator associated with \eqref{ibmnew3} is 
\begin{equation}\label{L1}
    \mathcal{L}_1^* f(x,t,v,m) = \frac{1}{T_m}\partial_m(g(m)f(x,t,v,m)) + \frac{\sigma^2}{2}\partial_m^2 f(x,t,v,m).
\end{equation}
As for the jump process \eqref{IBMmodel2}, the associated operator is
\begin{equation}\label{L2}
    \mathcal{L}_2^* f(x,t,v,m) = -\Lambda(m) f + \delta(m) \frac{1}{||\mathbb{V}_0||}\int_{\mathbb{V}_0} \int_{\mathbb{R}} \Lambda(m') f(t,x,v',m') dm' dv'.
\end{equation}
Using \eqref{L1} and \eqref{L2}, we derive the kinetic PDE for the individual based model \eqref{IBMmodel}--\eqref{IBMmodel2},
\begin{equation}\label{kineticeqn}
    \partial_t f + v \cdot \nabla_x f - \frac{1}{T_m} \partial_m [g(m)f] - \frac{\sigma^2}{2}\partial_m^2 f = -\Lambda(m) f + \delta(m) \frac{1}{||\mathbb{V}_0||}\int_{\mathbb{V}_0} \int_{\mathbb{R}} \Lambda(m') f(t,x,v',m') dm' dv',
\end{equation}
where $T_m$ is the adaptation time, $\delta(m)$ is the Dirac function, and $g(m)$, $\Lambda(m)$ are defined in \eqref{gm} and \eqref{Lambda} respectively.

We non-dimensionalize the kinetic equation \eqref{kineticeqn} by letting 
\begin{equation*}
    t = T_t \Tilde{t} , \quad x = L \Tilde{x}, \quad v = V_0 \Tilde{v}, \quad \Lambda = \Tilde{\Lambda}/T_{\Lambda},
\end{equation*}
where $T_t$ is the system time, $L$ is the characteristic spatial scale, and $V_0$ is the magnitude of the velocity, and $T_\Lambda$ is the characteristic tumbling time. 

Dropping the tilde sign for simplicity, the non-dimensionalized equation is
\begin{equation*}
    \frac{T_{\Lambda}}{T_t}\partial_{t} f + \frac{V_0T_{\Lambda}}{L}v \cdot \nabla_x f - \frac{T_{\Lambda}}{T_m} \partial_m [g(m)f] - \frac{T_{\Lambda}\sigma^2}{2} \partial_m^2 f = -\Lambda(m) f + \delta(m) \int_{\mathbb{V}} \int_{\mathbb{R}} \Lambda(m')  f dm' dv',
\end{equation*}
where $ v \in \mathbb{V}$. The velocity space $\mathbb{V}$ is $\{-1,1\}$ for $d=1$ and $\partial B(0,1)$ for $d \geq 2$ (the unit sphere) and $\, dv$ is the normalized surface measure.
Consider the kinetic equation \eqref{kineticeqn} under the scaling
\begin{equation}\label{eq:scaling}
    \epsilon = \frac{V_0T_{\Lambda}}{L},\quad \frac{T_{\Lambda}}{T_t} = \epsilon^{1+\mu}, \quad \frac{T_{\Lambda}}{T_m} = \epsilon^{\gamma},
\end{equation}
with $\epsilon$ small. Different $\mu$ and $\gamma$ correspond to different system time and adaptation time $T_m$. Now we illustrate the relations between the values of $\mu, \gamma$ and the numerical observations in Section 2.

\subsection{The scalings in the numerical simulations}

As in the numerical results detailed in Section \ref{sec:222}, we set $V_0=0.02$ and $\sigma=\sqrt{2}$. Since $\Lambda(m_t)$ satisfies \eqref{Lambda}, we choose $T_{\Lambda} = 1$. This leads to $\frac{T_\Lambda \sigma^2}{2} = 1$. The characteristic spatial scale is defined as the average spatial displacement during a certain time interval $(T_i, T_e)$. More precisely,
\begin{equation*}
    L(T_i, T_e) = \sqrt{\frac{1}{N} \sum_{n=1}^N \big|x^n(T_e)-x^n(T_i)\big|^2},
\end{equation*}
where the number of samples $N$ and position of the $n$th sample $x^n$ are defined the same as in Section \ref{sec:222}.

In Figure \ref{MSD}, the slopes of the MSD imply the dispersion patterns, which are quantified by the value of $\mu$. When particles transit from ballistic transport to normal diffusion, the value of $\mu$ changes with time accordingly. Therefore, in order to identify a particular dispersion pattern, one must consider a time interval with an almost constant MSD slope. To identity the value of $\mu$, let $T_{e2}>T_{e1}$ and $(T_i, T_{e1}) $, $(T_i, T_{e2}) $ be two intervals that exhibit the same dispersion pattern, i.e., they have the same MSD slope as in Figure \ref{MSD}. One can numerically calculate the characteristic lengths of these two intervals $L_1 = L(T_i, T_{e1})$, $L_2 = L(T_i, T_{e2})$ and use $T_{t1} = T_{e1}-T_i $, $T_{t2} = T_{e2}-T_i $ to determine $\mu$, while the value of $\epsilon$ may be different for these two intervals. Denote
\begin{equation*}
    \epsilon_j = \frac{V_0T_{\Lambda}}{L_j},\quad \frac{T_{\Lambda}}{T_{tj}} = \epsilon_j^{1+\mu}, 
\qquad j = 1, 2. 
\end{equation*}
From the above equations, given $L_j$ and $T_{tj}$ ($j=1,2$), one determines three unknowns from four equations. Note that the time interval is selected such that the MSD slope remains nearly constant, which ensures that the value $\mu$ stays approximately unchanged. Thus,
one can determine $\mu$ and $\epsilon_j$ by
\begin{equation*}
    \frac{T_{t2}}{T_{t1}} = \left(\frac{L_2}{L_1}\right)^{1+\mu} , \quad \epsilon_j = \left(\frac{T_{\Lambda}}{T_{tj}}\right)^{\frac{1}{1+\mu}}.
\end{equation*}
Then different $\epsilon$ correspond to different values of $\gamma$ such that
\begin{equation}\gamma_j = \log_{\epsilon_j} \frac{T_{\Lambda}}{T_m}.
\end{equation}

We consider two different time intervals: $[10^2, 10^3]$ and $[10^5, 10^6]$. In the interval $[10^2, 10^3]$, let $(T_i, T_{e1}) = (10^2, 6 \times 10^2)$ and $(T_i, T_{e2}) = (10^2, 10^3)$. Then, $T_{t1} = T_{e1} - T_i = 5 \times 10^2$ and $T_{t2} = T_{e2} - T_i = 9 \times 10^2$.
Similarly, in the interval $[10^5, 10^6]$, we set $(T_i, T_{e1}) = (10^5, 6 \times 10^5)$ and $(T_i, T_{e2}) = (10^5, 10^6)$. Thus, $T_{t1} = 5 \times 10^5$ and $T_{t2} = 9 \times 10^5$. For different $T_m$'s, we calculate the average spatial displacement $L(T_i, T_{ej})$ ($j=1,2$) and obtain the corresponding values of $1 + \mu$, $\epsilon_j$, and $\gamma_j$. These values are displayed in Tables \ref{table1} and \ref{table2}.

Tables \ref{table1} and \ref{table2} show that $\epsilon = \frac{V_0 T_{\Lambda}}{L}$ are small values ranging from $10^{-2}$ to $10^{-6}$. From Table \ref{table1}, in the interval $[10^2, 10^3]$, when $T_m = 10^2$, $10^3$, and $+\infty$, $\mu$ is approximately $0$, and particles exhibit ballistic transport at the population level. One can observe that both $\gamma_j$ ($j=1,2$) are above $1/2$ in these cases. On the other hand, when $T_m = 10^{-2}$ and $1$, $\mu$ is approximately $1$, indicating that particles perform normal diffusion. Moreover, both $\gamma_j$ ($j=1,2$) are below $0$. For the case when $T_m = 10$, the dispersion pattern can not be determined. Therefore, the range of $\gamma$ is related to the dispersion pattern of the particle population.
In Table \ref{table2}, for the longer time interval $[10^5, 10^6]$, similar observations can be made: when $\gamma$ is negative, particles perform normal diffusion, while when $\gamma_j > 1/2$ and $\mu \sim 0$, particles perform ballistic transportation.

\begin{table}[!htbp]
  \begin{center}
    \caption{The values of $L_{j}$, $\epsilon_j$, $1+\mu$ and $\gamma_j$ ($j=1,2$) in the interval $[10^2,10^3]$ for different $T_m$ when $T_{t1} = 5\times10^2$ and $T_{t2} = 9\times 10^2$. DP is the abbreviation for dispersion pattern, ND is for normal diffusion and BT is for ballistic transport.}
    \label{table1}
    \begin{tabular}{|c|c|c|c|c|c|c|c|c|} 
    \hline
      \textbf{$T_m$} & \textbf{$L_1$} & \textbf{$L_2$} & \textbf{$\epsilon_1$} & \textbf{$\epsilon_2$} & \textbf{$1+\mu$} & \textbf{$\gamma_1$} & \textbf{$\gamma_2$} & \text{DP}\\
      \hline
      $10^{-2}$ &$6.33\times10^{-1}$ & $8.63\times10^{-1}$ & $3.60\times10^{-2}$  & $2.77\times10^{-2}$ & $1.90$ & $-1.39$ & $-1.29$ & \text{ND}\\
      \hline
      $1$ & $7.18\times10^{-1}$ & $9.73\times10^{-1}$ & $4.06\times10^{-2}$  & $3.00\times10^{-2}$ & $1.94$ & $0$ & $0$ & \text{ND}\\
      \hline
      $10$ & $4.03$ & $5.81$ & $2.11\times10^{-2}$  & $1.47\times10^{-2}$ & $1.61$ & $0.60$ & $0.55$ & \text{--}\\
      \hline
      $10^2$ & $6.71$ & $1.16\times10^{1}$ & $3.10\times10^{-3}$ & $1.79\times10^{-3}$ & $1.08$ & $0.80$ & $0.73$ & \text{BT}\\
      \hline
      $10^3$ & $6.91$ & $1.21\times10^{1}$ & $2.70\times10^{-3}$ & $1.55\times10^{-3}$ & $1.05$ & $1.17$ & $1.07$ & \text{BT}\\
      \hline
      $+\infty$ & $6.87$ & $1.22\times10^{1}$  & $2.40\times10^{-3}$ & $1.34\times10^{-3}$ & $1.03$ & $+\infty$ & $+\infty$ & \text{BT}\\
      \hline
    \end{tabular}
  \end{center}
\end{table}
\begin{table}[!htbp]
  \begin{center}
    \caption{
    The values of $T_{tj}$, $L_{j}$, $\epsilon_j$, $1+\mu$ and $\gamma_j$ ($j=1,2$) in the interval $[10^5,10^6]$ for different $T_m$ when $T_{t1} = 5\times10^5$ and $T_{t2} = 9\times10^5$. Besides, DP is an abbreviation for dispersion pattern, ND is normal diffusion and BT is ballistic transport.}
    \label{table2}
    \begin{tabular}{|c|c|c|c|c|c|c|c|c|} 
    \hline
      \textbf{$T_m$} & \textbf{$L_1$} & \textbf{$L_2$}   & \textbf{$\epsilon_1$} & \textbf{$\epsilon_2$} & \textbf{$1+\mu$} & \textbf{$\gamma_1$} & \textbf{$\gamma_2$} & \text{DP} \\
      \hline
      $10^{-2}$ & $2.12\times10^{1}$ & $2.79\times10^{1}$ & $2.10\times10^{-3}$ & $1.59\times10^{-3}$ & $2.13$ & $-0.75$ & $-0.71$ & \text{ND}\\
      \hline
      $1$ & $2.3330\times10^{1}$ & $3.08\times10^{1}$ & $2.00\times10^{-3}$ & $1.51\times10^{-3}$ & $2.11$ & $0$ & $0$ & \text{ND}\\
      \hline
      $10$ & $1.57\times10^{2}$ & $2.06\times10^{2}$ & $2.40\times10^{-3}$ & $1.83\times10^{-3}$ & $2.17$ & $0.38$ & $0.37$ & \text{ND}\\
      \hline
      $10^2$ & $1.62\times10^{3}$ & $2.18\times10^{3}$ & $1.30\times10^{-3}$ & $1.00\times10^{-3}$ & $1.99$ & $0.70$ & $0.67$ & \text{-}\\
      \hline
      $10^3$ & $6.08\times10^{3}$ & $1.01\times10^{4}$ & $1.10\times10^{-5}$ & $6.57\times10^{-6}$ & $1.15$ & $0.60$ & $0.58$ & \text{BT}\\
      \hline
      $+\infty$ & $7.07\times10^{3}$ & $1.24\times10^{4}$ & $3.80\times10^{-6}$  & $2.17\times10^{-6}$ & $1.05$ & $+\infty$ & $+\infty$ & \text{BT}\\
      \hline
    \end{tabular}
  \end{center}
\end{table}

\subsection{The main theorem}
We consider the following rescaled PDE model  
\begin{equation}\label{pdemodel}
    \begin{cases}
    \epsilon^{1+\mu} \partial_t f_{\epsilon} + \epsilon v \cdot \nabla_x f_{\epsilon} - \epsilon^{\gamma} \partial_m [g(m)f_\epsilon] - \partial_m^2 f_{\epsilon} = -\Lambda(m) f_{\epsilon} + \delta(m) \int_{\mathbb{V}} \int_{\mathbb{R}} \Lambda(m') f_{\epsilon} dm' dv', \\
    f_{\epsilon}(x,t=0,v,m) = \delta(m) \rho_{in}(x).
\end{cases}
\end{equation}
The parameters $\mu$ and $\gamma$ satisfy $0 \leq \mu \leq 1$ and $\gamma \in \mathbb{R}$. Here $\gamma<0$ indicates that the adaptation is very fast or equivalently $T_m$ is small, while $\gamma$ being large indicates that the noise term dominates.

The main analytical result of this paper is as follows.
\begin{thm}\label{maintheorem}
    Let $f_{\epsilon}$ satisfy the equation \eqref{pdemodel} with the tumbling rate function $\Lambda$ given by \eqref{Lambda} and adaptation function $g(m)$ given by \eqref{gm}.
    \begin{itemize}
          \item [(a)] 
          When $\gamma > \frac{1}{2}$ and $\mu = 0$, there exists a function $\rho(x,t)$ such that
          \begin{equation*}
              \rho(x,t) = \lim_{\epsilon \to 0} \int_{\mathbb{V}} \int_{\mathbb{R}} f_{\epsilon}(x,t,v,m) \, dm \, dv
          \end{equation*}
          and $\rho(x,t)$ solves
          \begin{equation}\label{thmselfsimilar}
          \begin{cases}
              \rho(x,t) = 2\rho_{in} *\frac{1}{t^d}F(\frac{x}{t}) = \frac{2}{t^d} \int_{\mathbb{R}^d}\rho_{in}(y) F(\tfrac{x}{t} - y) dy, \\
              \rho(0,x) = \rho_{in}(x),
          \end{cases}
          \end{equation}
          where the fundamental solution $F(\frac{x}{t})$ is the inverse Fourier-Laplace transform of $H(is^{-1} \xi)$ defined in ~\eqref{laplacefourierF} and~\eqref{Hrho} and $*$ is the convolution symbol with respect to $x$. Furthermore, the associated \text{MSD} defined in \eqref{MSDgamma1/2} satisfies $\text{MSD} = C_0 t^2$, where $C_0$ is determined by the integral $C_0 = 2 \int_{\mathbb{R}^d} \int_{\mathbb{R}^d}|z|^2 \rho_{in}(y)F(z - y) dy dz$. In particular, for the one-dimensional case, $\rho(x,t)$ can be given explicitly by $\rho(x,t) = \int_{\mathbb{R}} \frac{\rho_{in}(x-y)}{\pi t\sqrt{1-(\frac{y}{t})^2}} \, dy$. 
          \item [(b)] 
          When $\gamma \leq 0$, we let $\mu = 1$. Then as $\epsilon \to 0$, there exists a function $\rho(x,t)$ such that
          \begin{equation*}
              \rho(x,t) = \lim_{\epsilon \to 0} \int_{\mathbb{V}} \int_{\mathbb{R}} f_{\epsilon}(t,x,m,v) \,dm \,dv,
          \end{equation*} 
          and $\rho(x,t)$ solves
          \begin{equation}\label{thmdiffusion}
              \begin{cases}
              \partial_t \rho - C \Delta_x \rho = 0, \\
              \rho(0,x) = \rho_{in}(x),
          \end{cases}
          \end{equation}
        where the constant $C$ can be computed explicitly and it depends on the dimension $d$.
       \end{itemize}
 \end{thm}

\begin{rmk}
The numerical results in Tables 1 and 2 are highly consistent with the analytical results, including those that are undetermined. In particular, in Table 1, for $T_m = 10^{-2}$ and $1$, the value of $\gamma$ is less than or equal to $0$ and $\mu$ is near 1, which, according to Theorem \ref{maintheorem} (b), corresponds to normal diffusion. Meanwhile, for $T_m = 10^2, 10^3, +\infty$, $\gamma$ is larger than $\frac{1}{2}$ and $\mu$ is close to 0,  indicating ballistic transport according to Theorem \ref{maintheorem} (a). The parameters $(\gamma, \mu)$ for the undetermined case $T_m = 10$ fall outside the ranges covered by Theorem~\ref{maintheorem}. Similar consistency is observed in Table 2.
\end{rmk}

\section{The proof of main results} \label{sec:proof}
In Section \ref{ODEsolving}, we start our analysis by applying the Laplace transform to the function $f_{\epsilon}$ with respect to the temporal variable $t$, followed by the Fourier transform with respect to the spatial variable $x$. To distinguish this transformed function from the original, $\hat{f_{\epsilon}}$ is the Fourier transform of $f_{\epsilon}$ with respect to $x$, and $\tilde{f_{\epsilon}}$ is the Laplace transform of $f_{\epsilon}$ with respect to $t$. This dual transformation transforms converts the partial differential equation \eqref{pdemodel} into the form presented in \eqref{fourierlaplace2}. Subsequently, we proceed to solve the resulting ordinary differential equation (ODE) for $\Tilde{\hat{f}}_{\epsilon}$ yielding the solutions $\Tilde{\hat{f}}_{\epsilon,\pm}$. These solutions are explicitly expressed in equations \eqref{expressionf+} and \eqref{expressionf-}. To obtain different limit equations for various ranges of $\gamma$, we calculate the asymptotic expansion of $\Tilde{\hat{f}}_{\epsilon,\pm}$ as $\gamma$ varies across different regimes.

For the limit equation when $\gamma > \frac{1}{2}$, we first derived the $\Tilde{\hat{\rho}}_{\epsilon}$ from $\hat{\rho}_{in}$ using the explicit solution \eqref{expressionf+} and \eqref{expressionf-}. We then derive the solution in a self-similar form, which arises naturally from the underlying assumptions and the specific properties of the system. Subsequently, we calculate the MSD of the particles, which manifests the ballistic transport of the particles on the macroscopic level.

When $\gamma \le 0$, to obtain the limit equation of \eqref{thmdiffusion} in Subsection \ref{derivationrho}, we use the explicit solutions \eqref{expressionf+} and \eqref{expressionf-} to derive the expression of $\Tilde{\hat{\rho}}_{\epsilon}$. Then, after performing the asymptotic expansion on the expression, we arrive at the limit equation of $\rho(x,t)$, which is a diffusion equation.

\subsection{The expression of $\Tilde{\hat{f}}_{\epsilon,\pm}$.}\label{ODEsolving}
Taking Fourier transform and Laplace transform of both sides of \eqref{pdemodel}, we have 
\begin{equation}\label{fourierlaplace2}
    \partial_m^2 \Tilde{\hat{f_{\epsilon}}} + \epsilon^{\gamma} \partial_m (g\tilde{\hat{f}}_{\epsilon})- [\epsilon^{1+\mu} s + i \epsilon \xi \cdot v + \Lambda(m)] \Tilde{\hat{f_{\epsilon}}} = -\epsilon^{1+\mu}\delta(m) \hat{\rho}_{in}(\xi) - \delta(m) \int_{\mathbb{V}} \int_{\mathbb{R}} \Lambda(m^{'}) \Tilde{\hat{f_{\epsilon}}} dm^{'} dv^{'}.
\end{equation}
Where $\Tilde{\hat{f}}_{\epsilon} = \Tilde{\hat{f}}_{\epsilon}(s,\xi,v,m)$. 
Defining $\Tilde{\hat{f}}_{\epsilon,+} = \Tilde{\hat{f}}|_{m \in (0,+\infty)}$ and $\Tilde{\hat{f}}_{\epsilon,-} = \Tilde{\hat{f}}|_{m \in (-\infty,0)}$, \eqref{fourierlaplace2} can be rewritten into the following two equations 
\begin{subequations}\label{Tildehatf2}
    \begin{equation}\label{m+2}
        \partial_m^2 \Tilde{\hat{f}}_{\epsilon,+} + \epsilon^{\gamma} \partial_m \Tilde{\hat{f}}_{\epsilon,+} - [\epsilon^{1+\mu} s + i \epsilon \xi \cdot v + 1] \Tilde{\hat{f}}_{\epsilon,+} = 0, 
    \end{equation}
    \begin{equation}\label{m-2}
        \partial_m^2 \Tilde{\hat{f}}_{\epsilon,-} -\epsilon^{\gamma} \partial_m \Tilde{\hat{f}}_{\epsilon,-} - [\epsilon^{1+\mu} s + i \epsilon \xi \cdot v] \Tilde{\hat{f}}_{\epsilon,-} = 0.
    \end{equation}
\end{subequations}

Integrating \eqref{fourierlaplace2} over the interval $[0-,0+]$ gives
\begin{equation}\label{derivate2}
     \partial_m \Tilde{\hat{f}}_{\epsilon,+}(0+) -  \partial_m \Tilde{\hat{f}}_{\epsilon,+}(0-) +  (\epsilon^{\gamma}\Tilde{\hat{f}}_{\epsilon,+}(0+)+ \epsilon^{\gamma} \Tilde{\hat{f}}_{\epsilon,+}(0-))  = - \epsilon^{1+\mu} \hat{\rho}_{in} - \int_{\mathbb{V}} \int_{\mathbb{R}}\Lambda(m^{'}) \Tilde{\hat{f_{\epsilon}}} dm^{'} dv^{'}.
\end{equation}
The continuity of $\Tilde{\hat{f}}_{\epsilon}$ at $m=0$ gives
\begin{equation}\label{continue2}
      \Tilde{\hat{f}}_{\epsilon,+}(0+) =   \Tilde{\hat{f}}_{\epsilon,+}(0-).
\end{equation}
Besides, it is required that the following limit be valid
\begin{equation}\label{limmtoinfty}
    \lim_{|m| \to \infty} \Tilde{\hat{f}}_{\epsilon} = 0.
\end{equation}

Considering the eigenvalue equation of \eqref{m+2}, we have
\begin{equation}\label{eigeqnlambda}
    \lambda^2 + \epsilon^{\gamma}\lambda - (\epsilon^{1+\mu} s + i \epsilon \xi \cdot v + 1) = 0.
\end{equation}
Letting $\lambda = a + i b$, we have 
\begin{numcases}{}
    a^2 - b^2 + \epsilon^{\gamma} a -(\epsilon^{1+\mu}s + 1) = 0, \label{a2} \\
    2ab + \epsilon^{\gamma} b -\epsilon \xi \cdot v = 0. \label{b2}
\end{numcases}
Similarly, considering the eigenvalue equation of \eqref{m-2}, we have 
\begin{equation}\label{eigeqnnu}
    \nu^2 - \epsilon^{\gamma} \nu -(\epsilon^{1+\mu}s +i \epsilon \xi \cdot v) = 0
\end{equation}
Letting $\nu = c + i d$, we get
\begin{numcases}{}
    c^2 - d^2 - \epsilon^{\gamma}c - \epsilon^{1+\mu}s = 0, \label{c2} \\
    2cd - \epsilon^{\gamma}d -\epsilon \xi \cdot v  = 0. \label{d2}
\end{numcases}

There are two solutions for the value with $\lambda$. The limitation \eqref{limmtoinfty} yields the one with the negative real part should be taken. As for the value of $\nu$, the one with the positive real part should be taken. Then, it follows that 
\begin{equation}\label{f+-2}
    \Tilde{\hat{f}}_{\epsilon,+} = C_+ e^{\lambda m},\quad 
    \Tilde{\hat{f}}_{\epsilon,-} = C_- e^{\nu m},
\end{equation}
where $\Re(\lambda) < 0$ and $\Re(\nu) > 0$.
substituting \eqref{f+-2} into \eqref{derivate2} and \eqref{continue2}, we have
\begin{equation*}
    C_+ = C_- = \frac{1}{\nu-\lambda - 2\epsilon^{\gamma}} \left(\epsilon^{1+\mu} \hat{\rho}_{in} + \ll \Tilde{\hat{f}}_{\epsilon,+} \gg\right),
\end{equation*}
where for simplicity of presentation, we introduce the notation
\begin{equation*}
    \ll \cdot \gg = \int_{\mathbb{V}} \int_{0}^{+\infty} \Lambda(m) \cdot  \,dm \, dv = \int_{\mathbb{V}} \int_{0}^{+\infty} \cdot \quad \,dm \, dv. 
\end{equation*}
Then,
\begin{equation}\label{expressionf+}
    \Tilde{\hat{f}}_{\epsilon,+} = \frac{e^{\lambda m}}{\nu-\lambda - 2\epsilon^{\gamma}} \left(\epsilon^{1+\mu} \hat{\rho}_{in} + \ll \Tilde{\hat{f}}_{\epsilon,+} \gg\right),
\end{equation}
\begin{equation}\label{expressionf-}
   \Tilde{\hat{f}}_{\epsilon,-} = \frac{e^{\nu m}}{\nu-\lambda - 2\epsilon^{\gamma}} \left(\epsilon^{1+\mu} \hat{\rho}_{in} + \ll \Tilde{\hat{f}}_{\epsilon,+} \gg\right).
\end{equation}

Integrating both sides of \eqref{expressionf+} in $m$ and $v$, we have
\begin{equation}\label{f+andrho}
    \ll \Tilde{\hat{f}}_{\epsilon,+} \gg = \epsilon^{1+\mu} \hat{\rho}_{in} \frac{\int_{\mathbb{V}} \frac{1}{-\lambda(\nu-\lambda-2\epsilon^{\gamma})} \, dv}{1-\int_{\mathbb{V}} \frac{1}{-\lambda(\nu-\lambda-2\epsilon^{\gamma})} \, dv}
\end{equation}
Substitute \eqref{f+andrho} into \eqref{expressionf+} and \eqref{expressionf-}, one has
\begin{equation}\label{f+rhoin}
    \Tilde{\hat{f}}_{\epsilon,+} = \frac{e^{\lambda m}}{\nu-\lambda-2\epsilon^{\gamma}} \frac{\epsilon^{1+\mu}\hat{\rho}_{in}}{1-\int_{\mathbb{V}} \frac{1}{-\lambda(\nu-\lambda-2\epsilon^{\gamma})} \, dv}.
\end{equation}
\begin{equation}\label{f-rhoin}
    \Tilde{\hat{f}}_{\epsilon,-} = \frac{e^{\nu m}}{\nu-\lambda-2\epsilon^{\gamma}} \frac{\epsilon^{1+\mu}\hat{\rho}_{in}}{1-\int_{\mathbb{V}} \frac{1}{-\lambda(\nu-\lambda-2\epsilon^{\gamma})} \, dv}.
\end{equation}

\paragraph{The expression of $\lambda = a + ib$}
\hfill \break
Substituting \eqref{b2} into \eqref{a2}, we obtain
\begin{equation*}
    [4(a^2 + \epsilon^{\gamma}a) + \epsilon^{2\gamma}](a^2 + \epsilon^{\gamma}a -\epsilon^{1+\mu}s-1) - \epsilon^2 |\xi \cdot v|^2 = 0.
\end{equation*}
Letting $A = a^2 + \epsilon^{\gamma}a$, we have
\begin{equation}\label{eqAgamma<3/2}
    4A^2 + (\epsilon^{2\gamma} - 4\epsilon^{1+\mu}s -4)A - \epsilon^{2\gamma}(\epsilon^{1+\mu}s+1) -\epsilon^2 |\xi \cdot v|^2 = 0.
\end{equation}
By calculating, we obtain 
\begin{align}\label{delta1gamma<3/2}
    \Delta_1 &={} (\epsilon^{2\gamma} - 4\epsilon^{1+\mu}s -4)^2 + 16\epsilon^{2\gamma}(\epsilon^{1+\mu}s+1) + 16\epsilon^2 |\xi \cdot v|^2
    \notag \\
    &={} 16 + 32 \epsilon^{1+\mu}s + 16 \epsilon^2 |\xi \cdot v|^2 + 16\epsilon^{2+2\mu}s^2 + 8\epsilon^{2\gamma}  + 8\epsilon^{2\gamma + 1+ \mu} s+\epsilon^{4\gamma}.
\end{align}
Solving equation \eqref{eqAgamma<3/2} yields
\begin{equation*}
    A = \frac{4\epsilon^{1+\mu}s + 4 -\epsilon^{2\gamma} \pm \sqrt{\Delta_1}}{8}.
\end{equation*}

Besides, solving the equation $A = a^2 + \epsilon^{\gamma}a$ yields \begin{equation*}
    a = \frac{-\epsilon^{\gamma} \pm \sqrt{\epsilon^{2\gamma} + 4A}}{2}.
\end{equation*}

From \eqref{b2} We can get the asymptotic expansion of $b$
\begin{equation}\label{expressionb}
    b =  \frac{\epsilon \xi \cdot v}{2a + \epsilon^{\gamma}} = \pm\frac{\epsilon \xi \cdot v}{\sqrt{\epsilon^{2\gamma} + 4A}}.
\end{equation}

\paragraph{The expression of $\nu = c + id$}
\hfill \break
Substituting \eqref{c2} into \eqref{d2} and denoting $B=c^2 - \epsilon^{\gamma}c$, we obtain
\begin{equation*}
    4B^2 + (\epsilon^{2\gamma} - 4\epsilon^{1+\mu}s)B -\epsilon^{2\gamma+1+\mu}s - \epsilon^2 |\xi \cdot v|^2 = 0.
\end{equation*}
$\Delta_2$ is 
\begin{align}\label{delta2gamma<3/2}
    \Delta_2 &={} (\epsilon^{2\gamma} - 4\epsilon^{1+\mu}s)^2 + 16\epsilon^{2\gamma+1+\mu}s + 16\epsilon^2 |\xi \cdot v|^2  \notag \\
    &={} 16 \epsilon^2 |\xi \cdot v|^2 + \epsilon^{4\gamma} + 8\epsilon^{2\gamma+1+\mu}s + 16 \epsilon^{2+2\mu}s^2.
\end{align}
By simple calculation, we have
\begin{equation}\label{B1/2gamma<3/2}
    B = \frac{4\epsilon^{1+\mu}s - \epsilon^{2\gamma} \pm \sqrt{\Delta_2}}{8}.
\end{equation}
Solving equation $c^2 - \epsilon^{\gamma}c - B =0$ yields
\begin{equation*}
    c = \frac{\epsilon^{\gamma} \pm \sqrt{\epsilon^{2\gamma}+4B}}{2}.
\end{equation*}
From \eqref{d2}, we obtain
\begin{equation}\label{expressiond}
    d = \frac{\epsilon \xi \cdot v}{2c - \epsilon^{\gamma}} = \pm \frac{\epsilon \xi \cdot v}{\sqrt{\epsilon^{2\gamma}+4B}}.
\end{equation}

\subsection{The derivation of the limit equation when $\gamma > \frac{1}{2}$.}
We first derived the limit equation for $\rho(x,t)$, and then, we analyzed the one-dimensional case, from which we were able to directly write out the explicit expression for $\rho(x,t)$.
\subsubsection{The derivation of the limit equation}
When $\gamma > \frac{1}{2}$, we set $\mu = 0$. From \eqref{f+rhoin} and \eqref{f-rhoin}, we have
\begin{equation}\label{transportrho}
    \Tilde{\hat{\rho}}_{\epsilon} = \int_{\mathbb{V}} \int_{\mathbb{R}} \Tilde{\hat{f}}_{\epsilon} \, dm \, dv =\int_{\mathbb{V}} \left(\int_0^{+\infty} \Tilde{\hat{f}}_{\epsilon,+} \, dm + \int_{-\infty}^0 \Tilde{\hat{f}}_{\epsilon,-} \, dm\right) \, dv =\epsilon \hat{\rho}_{in} \frac{\int_{\mathbb{V}} \frac{1}{-\lambda \nu}\left(1+\frac{2\epsilon^{\gamma}}{\nu-\lambda-2\epsilon^{\gamma}}\right) \, dv}{1-\int_{\mathbb{V}} \frac{1}{-\lambda\left(\nu-\lambda-2\epsilon^{\gamma}\right)} \, dv}.
\end{equation}

Next, we need to calculate the asymptotic expansion of the RHS of \eqref{transportrho}. 
The asymptotic expansions of $\lambda$ and $\nu$ are respectively
\begin{equation*}
    \lambda = -1-\epsilon \frac{s}{2}-\epsilon^{\gamma}\frac{1}{2} + O(\epsilon^2) - i\left[\epsilon \frac{\xi \cdot v}{2}+O(\epsilon^2)\right].
\end{equation*}
\begin{equation*}
    \nu = \epsilon^{\frac{1}{2}} \frac{\sqrt{2}}{2} \left(s+\sqrt{s^2+|\xi \cdot v|^2}\right)^{\frac{1}{2}} + \frac{\epsilon^{\gamma}}{2} + O(\epsilon^{2\gamma-\frac{1}{2}}) + i\left[\epsilon^{\frac{1}{2}} \frac{\frac{\sqrt{2}}{2} \xi \cdot v}{\left(s+\sqrt{s^2+|\xi \cdot v|^2}\right)^{\frac{1}{2}}} + O(\epsilon^{2\gamma-\frac{1}{2}})\right].
\end{equation*}

\paragraph{The asymptotic expansion of $\int_{\mathbb{V}} \frac{1}{\lambda \nu}\left(1+\frac{2\epsilon^{\gamma}}{\nu-\lambda-2\epsilon^{\gamma}}\right) \, dv$.}
By simple calculation, we have
\begin{align*}
    \nu -\lambda -2\epsilon^{\gamma} &={} 1+\epsilon^{\frac{1}{2}} \left(s+\sqrt{s^2+|\xi \cdot v|^2}\right)^{\frac{1}{2}} - \epsilon^{\gamma} + O(\epsilon^{2\gamma - \frac{1}{2}}) \notag \\
    &+{} i\left[\epsilon^{\frac{1}{2}} \frac{\frac{\sqrt{2}}{2}\xi \cdot v}{\left(s+\sqrt{s^2+|\xi \cdot v|^2}\right)^{\frac{1}{2}}} + O(\epsilon^{2\gamma-\frac{1}{2}})\right].
\end{align*}
The reciprocal of $\nu-\lambda-2\epsilon^{\gamma}$ is 
\begin{equation*}
    \frac{1}{\nu-\lambda-2\epsilon^{\gamma}} = 1+O(\epsilon^{\frac{1}{2}}) + i O(\epsilon^{\frac{1}{2}}).
\end{equation*}
Then,
\begin{equation*}
    1+\frac{2\epsilon^{\gamma}}{\nu-\lambda-2\epsilon^{\gamma}} = 1 + O(\epsilon^{\gamma}) + iO(\epsilon^{\gamma+\frac{1}{2}}).
\end{equation*}
As a result, we have
\begin{equation*}
    \frac{1}{-\lambda \nu}\left(1+\frac{2\epsilon^{\gamma}}{\nu-\lambda-2\epsilon^{\gamma}}\right)  \sim \frac{1}{-\lambda \nu}.
\end{equation*}
Multiplying $-\lambda$ and $\nu$, one gets
\begin{equation*}
    -\lambda \nu = \epsilon^{\frac{1}{2}} \frac{\sqrt{2}}{2} \left(s+\sqrt{s^2+|\xi \cdot v|^2}\right)^{\frac{1}{2}} + \epsilon^{\gamma} \frac{1}{2} + O(\epsilon^{2\gamma-\frac{1}{2}}) + i\left[ \epsilon^{\frac{1}{2}} \frac{\frac{\sqrt{2}}{2} \xi \cdot v}{\left(s+\sqrt{s^2+|\xi \cdot v|^2}\right)^{\frac{1}{2}}} + O(\epsilon^{2\gamma-\frac{1}{2}})\right].
\end{equation*}
The square of the magnitude of $-\lambda \nu$ is
\begin{equation*}
    \Re^2(-\lambda \nu) + \Im^2 (-\lambda \nu) = \epsilon \left[\frac{1}{2}\left(s+\sqrt{s^2+|\xi \cdot v|^2}\right) + \frac{1}{2}\frac{|\xi \cdot v|^2}{s+\sqrt{s^2 + |\xi \cdot v|^2}}\right] + O(\epsilon^{\gamma+\frac{1}{2}}). 
\end{equation*}
Then we can get the reciprocal 
\begin{equation*}
    \left[\Re^2(-\lambda \nu) + \Im^2(-\lambda \nu)\right]^{-1} =  \epsilon^{-1} \left[\frac{1}{2}\left(s+\sqrt{s^2+|\xi \cdot v|^2}\right) + \frac{1}{2}\frac{|\xi \cdot v|^2}{s+\sqrt{s^2 + |\xi \cdot v|^2}}\right]^{-1} + O(\epsilon^{\gamma-\frac{3}{2}}). 
\end{equation*}
As a consequence, the reciprocal of $-\lambda \nu$ is 
\begin{align*}
    \frac{1}{-\lambda \nu} &={} \epsilon^{-\frac{1}{2}} \frac{\sqrt{2}\left(s+\sqrt{s^2+|\xi \cdot v|^2}\right)^{\frac{1}{2}}}{\left[\left(s+\sqrt{s^2+|\xi \cdot v|^2}\right) + \frac{|\xi \cdot v|^2}{s+\sqrt{s^2 + |\xi \cdot v|^2}}\right]} + O(\epsilon^{\gamma-1}) \notag \\
    &-{} i\left[\epsilon^{-\frac{1}{2}} \frac{\frac{\sqrt{2}\xi \cdot v}{\left(s+\sqrt{s^2+|\xi \cdot v|^2}\right)^{\frac{1}{2}}}}{\left[\left(s+\sqrt{s^2+|\xi \cdot v|^2}\right) + \frac{|\xi \cdot v|^2}{s+\sqrt{s^2 + |\xi \cdot v|^2}}\right]} + O(\epsilon^{2\gamma-\frac{3}{2}})\right] \notag \\
    &={} \epsilon^{-\frac{1}{2}} \frac{\sqrt{2}\left(s+\sqrt{s^2+|\xi \cdot v|^2} \right)^{\frac{3}{2}}}{\left(s+\sqrt{s^2+|\xi \cdot v|^2}\right)^2 + |\xi \cdot v|^2} + O(\epsilon^{\gamma-1}) \notag \\
    &-{} i\left[\epsilon^{-\frac{1}{2}} \frac{\sqrt{2}\xi \cdot v\left(s+\sqrt{s^2 + |\xi \cdot v|^2}\right)^{\frac{1}{2}}}{\left(s+\sqrt{s^2 + |\xi \cdot v|^2}\right)^2 + |\xi \cdot v|^2} + O(\epsilon^{2\gamma-\frac{3}{2}}) \right].
\end{align*}
The integral of $\frac{1}{-\lambda \nu}$ is 
\begin{equation}\label{numeratorgamma1}
    \int_{\mathbb{V}} \frac{1}{-\lambda \nu} \, dv = \epsilon^{-\frac{1}{2}} \sqrt{2} \int_{\mathbb{V}} \frac{\left(s+\sqrt{s^2+|\xi \cdot v|^2}\right)^{\frac{3}{2}}}{\left(s+\sqrt{s^2+|\xi \cdot v|^2}\right)^2 + |\xi \cdot v|^2} + O(\epsilon^{\gamma-1}) + iO(\epsilon^{2\gamma-\frac{3}{2}}).
\end{equation}

\paragraph{The asymptotic expansion of $1-\int_{\mathbb{V}} \frac{1}{-\lambda\left(\nu-\lambda-2\epsilon^{\gamma}\right)} \, dv$.}
By simple computation, we get 
\begin{equation*}
    -\lambda(\nu-\lambda-2\epsilon^{\gamma}) = 1+\epsilon^{\frac{1}{2}} \frac{\sqrt{2}}{2} \left(s+\sqrt{s^2+|\xi \cdot v|^2}\right)^{\frac{1}{2}} + O(\epsilon^{\gamma}) + i\left[ \epsilon^{\frac{1}{2}} \frac{\frac{\sqrt{2}}{2} \xi \cdot v}{\left(s+\sqrt{s^2+|\xi \cdot v|^2}\right)^{\frac{1}{2}}} + O(\epsilon^{2\gamma-\frac{1}{2}})\right].
\end{equation*}
\begin{equation*}
    \Re^2[-\lambda(\nu-\lambda-2\epsilon^{\gamma})] + \Im^2[-\lambda(\nu-\lambda-2\epsilon^{\gamma})] = 1+\epsilon^{\frac{1}{2}} \sqrt{2} \left(s+\sqrt{s^2+|\xi \cdot v|^2}\right)^{\frac{1}{2}} + O(\epsilon).
\end{equation*}
Furthermore, we have
\begin{equation*}
    \{\Re^2[-\lambda(\nu-\lambda-2\epsilon^{\gamma})]+\Im^2[-\lambda(\nu-\lambda-2\epsilon^{\gamma})]\}^{-1} = 1-\epsilon^{\frac{1}{2}} \sqrt{2} \left(s+\sqrt{s^2+|\xi \cdot v|^2}\right)^{\frac{1}{2}} + O(\epsilon).
\end{equation*}
Additionally, we get the reciprocal of $-\lambda(\nu-\lambda-2\epsilon^{\gamma})$ writing
\begin{equation*}
    \frac{1}{-\lambda(\nu-\lambda-2\epsilon^{\gamma})} = 1 - \epsilon^{\frac{1}{2}} \sqrt{2} \left(s+\sqrt{s^2+|\xi \cdot v|^2}\right)^{\frac{1}{2}} + O(\epsilon^{\gamma}) - i\left[\epsilon^{\frac{1}{2}} \frac{\frac{\sqrt{2}}{2} \xi \cdot v}{\left(s+\sqrt{s^2+|\xi \cdot v|^2}\right)^{\frac{1}{2}}} + O(\epsilon^{2\gamma-\frac{1}{2}})\right].
\end{equation*}

As a result, 
\begin{equation}\label{denominatorgamma1}
    1-\int_{\mathbb{V}} \frac{\, dv}{-\lambda(\nu-\lambda-2\epsilon^{\gamma})} = \epsilon^{\frac{1}{2}} \frac{\sqrt{2}}{2} \int_{\mathbb{V}} \left(s+\sqrt{s^2+|\xi \cdot v|^2}\right)^{\frac{1}{2}} \, dv + O(\epsilon^{\gamma}) + iO(\epsilon^{2\gamma-\frac{1}{2}}).
\end{equation}

\paragraph{The expression of $\Tilde{\hat{\rho}}_{\epsilon}$.} 
From \eqref{denominatorgamma1}, one can get
\begin{equation}\label{denominatorgamma1new}
    \left[1-\int_{\mathbb{V}} \frac{\, dv}{-\lambda(\nu-\lambda-2\epsilon^{\gamma})}\right]^{-1} = \epsilon^{-\frac{1}{2}} \sqrt{2} \left[\int_{\mathbb{V}} \left(s+\sqrt{s^2+|\xi \cdot v|^2}\right) \, dv\right]^{-1} + O(\epsilon^{\gamma-1}) + iO(\epsilon^{2\gamma-\frac{3}{2}}).
\end{equation}
Then, by \eqref{numeratorgamma1} and \eqref{denominatorgamma1new}, we have
\begin{equation*}
    \frac{\int_{\mathbb{V}} \frac{1}{-\lambda \nu} \, dv}{1-\int_{\mathbb{V}} \frac{\, dv}{-\lambda(\nu-\lambda-2\epsilon^{\gamma})}} =2 \epsilon^{-1} \frac{ \displaystyle  \int_{\mathbb{V}} \frac{\left(s+\sqrt{s^2+|\xi \cdot v|^2}\right)^{\frac{3}{2}}}{\left(s+\sqrt{s^2+|\xi \cdot v|^2}\right)^2 + |\xi \cdot v|^2} \, dv}{\displaystyle \int_{\mathbb{V}} \left(s+\sqrt{s^2+|\xi \cdot v|^2}\right)^{\frac{1}{2}} \, dv} + O(\epsilon^{\gamma-\frac{3}{2}}) + iO(\epsilon^{2\gamma-2}).
\end{equation*}
Moreover, we can get the asymptotic expansion of \eqref{transportrho}
\begin{align}\label{transportasy}
    \Tilde{\hat{\rho}}_{\epsilon} = \epsilon \hat{\rho}_{in} \frac{\int_{\mathbb{V}} \frac{1}{-\lambda \nu} \, dv}{1-\int_{\mathbb{V}} \frac{\, dv}{-\lambda(\nu-\lambda-2\epsilon^{\gamma})}}  &={} 2\hat{\rho}_{in} \frac{\displaystyle \int_{\mathbb{V}} \frac{\left(s+\sqrt{s^2+|\xi \cdot v|^2}\right)^{\frac{3}{2}}}{\left(s+\sqrt{s^2+|\xi \cdot v|^2}\right)^2 + |\xi \cdot v|^2} \, dv}{\displaystyle \int_{\mathbb{V}} \left(s+\sqrt{s^2+|\xi \cdot v|^2}\right)^{\frac{1}{2}} \, dv} + O(\epsilon^{\gamma-\frac{1}{2}}) + iO(\epsilon^{2\gamma-1}).
\end{align}

Letting $\Tilde{\hat{\rho}} = \lim_{\epsilon \to 0} \Tilde{\hat{\rho}}_{\epsilon}$. As $\epsilon$ to $0$, we can get the limit equation of $\Tilde{\hat{\rho}}$
\begin{equation*}
    \Tilde{\hat{\rho}} = 2\hat{\rho}_{in} \frac{\displaystyle \int_{\mathbb{V}} \frac{\left(s+\sqrt{s^2+|\xi \cdot v|^2}\right)^{\frac{3}{2}}}{\left(s+\sqrt{s^2+|\xi \cdot v|^2}\right)^2 + |\xi \cdot v|^2} \, dv}{\displaystyle \int_{\mathbb{V}} \left(s+\sqrt{s^2+|\xi \cdot v|^2}\right)^{\frac{1}{2}} \, dv}.
\end{equation*}

For a function $\frac{1}{t^d}F(\frac{x}{t})$, $x \in \mathbb{R}^d$, let $y = \frac{x}{t}$, after undergoing Fourier and Laplace transformations, one can obtain
\begin{align}\label{laplacefourierF}
    \mathscr{F}\left[\mathscr{L}[\frac{1}{t^d}F(\frac{x}{t})]\right] &={} \int_{\mathbb{R}^d} \int_0^{\infty} e^{-i\xi \cdot x - st} \frac{1}{t^d}F(\frac{x}{t}) \, dt \, dx = 
    \int_{\mathbb{R}^d} \int_0^{\infty} e^{-(i \xi \cdot \frac{x}{t}+s)t} F(\frac{x}{t}) \, dt \, d(\frac{x}{t}) \notag \\
    &={} \int_{\mathbb{R}^d} \int_0^{\infty} e^{-(i \xi \cdot y + s)t} F(y) \, dt \, dy = \int_{\mathbb{R}^d} \frac{F(y)}{i \xi \cdot y + s} \, dy
    = s^{-1} \int_{\mathbb{R^d}} \frac{F(y)}{1+\frac{i \xi}{s} \cdot y} \, dy \notag \\ 
    &\triangleq{} s^{-1} H(i s^{-1} \xi),
\end{align}

Moreover, from \eqref{laplacefourierF} the limit equation of $\Tilde{\hat{\rho}}$ can be represented 
\begin{align}
    \Tilde{\hat{\rho}} &={} 2\hat{\rho}_{in} \frac{\displaystyle \int_{\mathbb{V}} \frac{s^{\frac{3}{2}}\left(1+\sqrt{1+s^{-2}|\xi \cdot v|^2}\right)^{\frac{3}{2}}}{s^2\left(1+\sqrt{1+s^{-2}|\xi \cdot v|^2}\right)^2 + s^{-2}|\xi \cdot v|^2} \, dv}{\displaystyle \int_{\mathbb{V}} s^{\frac{1}{2}}\left(1+\sqrt{1 + s^{-2}|\xi \cdot v|^2}\right)^{\frac{1}{2}} \, dv} \notag \\
    &={} 2\hat{\rho}_{in} s^{-1} \frac{\displaystyle \int_{\mathbb{V}} \frac{\left(1+\sqrt{1+(s^{-1}|\xi|)^2 |e_{\xi} \cdot v|^2}\right)^{\frac{3}{2}}}{\left(1+\sqrt{1+(s^{-1}|\xi|)^2 |e_{\xi} \cdot v|^2}\right)^2 + (s^{-1}|\xi|)^2 |e_{\xi} \cdot v|^2} \, dv}{\displaystyle \int_{\mathbb{V}} \left(1+\sqrt{1 + (s^{-1}|\xi|)^2 |e_{\xi} \cdot v|^2}\right)^{\frac{1}{2}} \, dv} \notag \\
    &={}2\hat{\rho}_{in} s^{-1} \frac{\displaystyle \int_{\mathbb{V}} \frac{\left(1+\sqrt{1+(s^{-1}\xi) \cdot (s^{-1} \xi) |e_{\xi} \cdot v|^2}\right)^{\frac{3}{2}}}{\left(1+\sqrt{1+(s^{-1}\xi) \cdot (s^{-1} \xi) |e_{\xi} \cdot v|^2}\right)^2 + (s^{-1}\xi) \cdot (s^{-1} \xi) |e_{\xi} \cdot v|^2} \, dv}{\displaystyle \int_{\mathbb{V}} \left(1+\sqrt{1 + (s^{-1}\xi) \cdot (s^{-1} \xi) |e_{\xi} \cdot v|^2}\right)^{\frac{1}{2}} \, dv} \\
    &={} 2\hat{\rho}_{in} s^{-1} \frac{\displaystyle \int_{\mathbb{V}} \frac{\left(1+\sqrt{1-\eta \cdot \eta |e_{\xi} \cdot v|^2}\right)^{\frac{3}{2}}}{\left(1+\sqrt{1-\eta \cdot \eta |e_{\xi} \cdot v|^2}\right)^2 - \eta \cdot \eta |e_{\xi} \cdot v|^2} \, dv}{\displaystyle \int_{\mathbb{V}} \left(1+\sqrt{1 - \eta \cdot \eta |e_{\xi} \cdot v|^2}\right)^{\frac{1}{2}} \, dv} \notag \\
    &\triangleq{} 2\hat{\rho}_{in} s^{-1} H(\eta). \label{Hrho}
\end{align}
where $\eta = is^{-1}\xi$ and $\xi = |\xi|e_{\xi}$.

Consequently, there exists a scaling function $F(\frac{x}{t})$ such that, upon performing the inverse Laplace transform followed by the inverse Fourier transform on $s^{-1}H$, the result is equivalent to $\frac{1}{t^d}F(\frac{x}{t})$. It follows that the limit equation for $\rho(x,t)$ is 
\begin{equation*}
    \rho(x,t) = 2\rho_{in}*\frac{1}{t^d}F(\frac{x}{t}),
\end{equation*}
where $*$ is the convolution symbol.

The functional form of $F(\frac{x}{t})$ suggests a scaling variable $\frac{x}{t}$, meaning that a characteristic spatial scale on which the density changes $x$ scales with time as $x \propto t$, causing the MSD to be of the form $t^2$, which is indicative ballistic transport on the macroscopic scale. To verify this idea, we calculate the MSD expression.

\paragraph{The calculation of the MSD}
\begin{equation*}
    {\rm MSD} = \int_{\mathbb{R}^d} |x|^2 \rho(x, t) dx
    = 2 \int_{\mathbb{R}^d}
    |x|^2 \rho_{in}*\frac{1}{t^d}F(\frac{x}{t}) dx
    = \frac{2}{t^d}
    \int_{\mathbb{R}^d} \int_{\mathbb{R}^d}
    |x|^2 \rho_{in}(y)
    F(x/t - y) dy dx.
\end{equation*}
Make the change of variables $z = x/t$ for $t > 0$. Then
\begin{align}\label{MSDgamma1/2}
\int_{\mathbb{R}^d} |x|^2 \rho(x, t) dx
= 2 t^{2} \int_{\mathbb{R}^d} \int_{\mathbb{R}^d}
   |z|^2 \rho_{in}(y)
   F(z - y) dy dz
= C_0 t^{2}. 
\end{align}
where $C_0 = 2 \int_{\mathbb{R}^d} \int_{\mathbb{R}^d}|z|^2 \rho_{in}(y)F(z - y) dy dz$.


\subsubsection{The one-dimensional case for $\gamma>\frac{1}{2}$}
\paragraph{Theoretical analysis}
Considering the case where the spatial coordinate $x$ and the velocity $v$ are confined to one dimension, we have $\xi \in \mathbb{R}$ and $v \in \{ \pm 1\}$.
\begin{align*}
    \Tilde{\hat{\rho}} &={} 2\hat{\rho}_{in} \frac{\displaystyle \int_{\mathbb{V}} \frac{\left(s+\sqrt{s^2+|\xi v|^2}\right)^{\frac{3}{2}}}{\left(s+\sqrt{s^2+|\xi v|^2}\right)^2 + |\xi v|^2} \, dv}{\displaystyle \int_{\mathbb{V}} \left(s+\sqrt{s^2+|\xi v|^2}\right)^{\frac{1}{2}} \, dv}
    = 2 \hat{\rho}_{in} \frac{ \frac{\left(s+\sqrt{s^2+\xi^2}\right)^{\frac{3}{2}}}{\left(s+\sqrt{s^2+\xi^2}\right)^2 + \xi^2}}{ \left(s+\sqrt{s^2+\xi^2}\right)^{\frac{1}{2}}} \notag \\
    &={} 2 \hat{\rho}_{in} \frac{1}{s} \frac{1+\sqrt{1+ (s^{-1}\xi)^2}}{\left(1+\sqrt{1+(s^{-1}\xi)^2}\right)^2 + (s^{-1}\xi)^2} = 2 \hat{\rho}_{in} \frac{1}{s} \frac{1+\sqrt{1-(i s^{-1}\xi)^2}}{\left(1+\sqrt{1-(i s^{-1}\xi)^2}\right)^2 - (i s^{-1}\xi)^2} \notag \\
    &\triangleq{} 2\hat{\rho}_{in} \frac{1}{s} H(is^{-1} \xi).
\end{align*}

For the function $H(is^{-1} \xi)$, using the Sokhotsky-Weierstrass theorem, one can obtain the corresponding function $F(\frac{x}{t})$ without performing the inverse Laplace and Fourier transforms \cites{levywalks, PhysRevE.91.022131}. It's noted that
\begin{equation*}
    F(\frac{x}{t}) = -\frac{1}{\pi \frac{x}{t}} \lim_{\delta \to 0} \Im [H(-\frac{1}{\frac{x}{t}+i\delta})] = -\frac{1}{\pi y} \lim_{\delta \to 0} \Im [H(-\frac{1}{y+i\delta})],
\end{equation*}
where $y = \frac{x}{t}$.

Then, we have
\begin{equation*}
    H(-\frac{1}{y+i\delta}) = \frac{1}{s} \frac{1+\sqrt{1-(y+i\delta)^{-2}}}{\left(1+\sqrt{1-(y+i\delta)^{-2}}\right)^2 -(y+i\delta)^{-2}}.
\end{equation*}
By simple computation,
\begin{equation*}
    (y+i\delta)^{-2} = \frac{1}{y^2-\delta^2+i2y\delta} = \frac{y^2-\delta^2}{(y^2+\delta^2)^2} - i\frac{2y\delta}{(y^2+\delta^2)^2}.
\end{equation*}
\begin{equation*}
    1-(y+i\delta)^{-2} = 1-\frac{(y^2-\delta^2)}{(y^2+\delta^2)^2} + i \frac{2y\delta}{(y^2+\delta^2)^2}.
\end{equation*}
Letting
\begin{equation*}
    z_0 = a_0 + ib_0 = \sqrt{1-(y+i\delta)^{-2}}.
\end{equation*}
We have
\begin{equation*}
    (y+i\delta)^{-2} = 1-a_0^2+b_0^2-i2a_0b_0.
\end{equation*}
It can be obtained 
\begin{equation}\label{z0mol}
    r_0 = |1-(y+i\delta)^{-2}| =  \sqrt{\left(1-\frac{(y^2-\delta^2)}{(y^2+\delta^2)^2}\right)^2 + \frac{4y^2\delta^2}{(y^2+\delta^2)^4}}.
\end{equation}
\begin{equation}\label{a0}
    a_0 = \pm \sqrt{\frac{|r_0|+1-\frac{(y^2-\delta^2)}{(y^2+\delta^2)^2}}{2}}.
\end{equation}
\begin{equation}\label{b0}
    b_0 = \pm \frac{y}{|y|} \sqrt{\frac{|r_0|-1+\frac{(y^2-\delta^2)}{(y^2+\delta^2)^2}}{2}}.
\end{equation}

Therefore, $H(-\frac{1}{y+i\delta})$ can be expressed as
\begin{align*}
    H(-\frac{1}{y+i\delta}) &={} \frac{1+a_0+ib_0}{\left(1+a_0+ib_0\right)^2 - (1-a_0^2+b_0^2-i2a_0b_0)}\notag\\
    &={} \frac{1+a_0+ib_0}{2a_0^2-2b_0^2+2a_0 + i(2b_0 + 4a_0b_0)} \notag \\
    &={} \frac{(1+a_0+ib_0)[2a_0^2-2b_0^2+2a_0 - i(2b_0 + 4a_0b_0)]}{(2a_0^2-2b_0^2+2a_0)^2 + (2b_0 + 4a_0b_0)^2}.
\end{align*}
Then, $\Im[H(-\frac{1}{y+i\delta})]$ is 
\begin{align}\label{ImH}
    \Im[H(-\frac{1}{y+i\delta})] &={}   \frac{-(1+a_0)(2b_0 + 4a_0b_0) + b_0(2a_0^2-2b_0^2+2a_0)}{(2a_0^2-2b_0^2+2a_0)^2 + (2b_0 + 4a_0b_0)^2} \notag \\
    &={} \frac{-2b_0(a_0^2+b_0^2+2a_0+1)}{(2a_0^2-2b_0^2+2a_0)^2 + 4b_0^2(1 + 2a_0)^2}.
\end{align}

When $\delta \to 0$, from \eqref{z0mol}, we obtain
\begin{equation*}
    r \triangleq \lim_{\delta \to 0}|r_0| = |1-y^{-2}|,
\end{equation*}
which means
\begin{equation*}
    r = \begin{cases}
                1-y^{-2}, \quad y^2 \geq 1, \\
                y^{-2} -1, \quad y^2 < 1.
           \end{cases}
\end{equation*}
Then, by \eqref{a0} and \eqref{b0} one can get
\begin{equation*}
    \alpha \triangleq \lim_{\delta \to 0} a_0 = \begin{cases}
                    \pm \sqrt{1-y^{-2}}, \quad y^2 \geq 1, \\
                    0, \quad y^2 < 1.
           \end{cases}
\end{equation*}
and 
\begin{equation*}
    \beta \triangleq \lim_{\delta \to 0} b_0 = \begin{cases}
                    0, \quad y^2 \geq 1, \\
                    \pm \frac{y}{|y|} \sqrt{y^{-2}-1}, \quad y^2 < 1.
           \end{cases}
\end{equation*}
Furthermore, when $y^2 > 1$, we find that $\lim_{\delta \to 0} b_0 = 0$ and $\lim_{\delta \to 0} b_1 = 0$. Consequently, this leads to the imaginary part of $H(-\frac{1}{y+i\delta})$ vanishing as $\delta \to 0$. To avoid this situation, we assume that $y^2 < 1$ is satisfied. Then, we have
\begin{equation*}
    r = y^{-2} - 1.
\end{equation*}
\begin{equation*}
    \alpha = 0.
\end{equation*}
\begin{equation*}
    \beta = \pm \frac{y}{|y|} \sqrt{y^{-2} - 1}.
\end{equation*}

As a result, from \eqref{ImH}, $\lim_{\delta \to 0} \Im[H(-\frac{1}{y+i\delta})]$ becomes
\begin{equation*}
     \lim_{\delta \to 0}\Im[H(-\frac{1}{y+i\delta})] = \frac{-2\beta(\beta^2+1)}{4\beta^4 + 4\beta^2} = -\frac{1}{2\beta}.
\end{equation*}

We can obtain different results using the formula calculated above based on the sign of $\beta_0$. We will discuss the cases separately
\paragraph{Case 1: $\beta_0 = \frac{y}{|y|} \sqrt{y^{-2}-1}$}.
By computation, we have
\begin{equation*}
     \lim_{\delta \to 0}\Im[H(-\frac{1}{y+i\delta})] = -\frac{1}{2\beta_0} = -\frac{1}{2\frac{y}{|y|}\sqrt{y^{-2}-1}}= -\frac{|y|^2}{2y\sqrt{1-y^2}} = - \frac{y}{2\sqrt{1-y^2}}.
\end{equation*}
By Sokhotsky-Weierstrass theorem, we obtain
\begin{equation*}
    F(\frac{x}{t}) = -\frac{1}{\pi y} \lim_{\delta \to 0}\Im[H(-\frac{1}{y+i\delta})] = \frac{1}{2\pi \sqrt{1-y^2}}.
\end{equation*}
Then, the limit equation for $\rho$ is
\begin{equation*}
    \rho(x,t) = 2\rho_{in}*\frac{1}{t}F(\frac{x}{t}) = \frac{1}{\pi t}\rho_{in}(x)*\frac{1}{\sqrt{1-(\frac{x}{t})^2}}.
\end{equation*}

\paragraph{Case 2: $\beta_0 = -\frac{y}{|y|} \sqrt{y^{-2}-1}$}.

Through calculations similar to those in Case 1, we can obtain
\begin{equation*}
     \lim_{\delta \to 0}\Im[H(-\frac{1}{y+i\delta})] = -\frac{1}{2\beta_0} = \frac{y}{2\sqrt{1-y^2}}.
\end{equation*}
By Sokhotsky-Weierstrass theorem, we obtain
\begin{equation*}
    F(\frac{x}{t}) = -\frac{1}{\pi y} \lim_{\delta \to 0}\Im[H(-\frac{1}{y+i\delta})] =-\frac{1}{2 \pi \sqrt{1-y^2}}.
\end{equation*}
Then, the limit equation for $\rho$ is
\begin{equation*}
    \rho(x,t) = 2\rho_{in}*\frac{1}{t}F(\frac{x}{t}) = -\frac{1}{2 \pi t}\rho_{in}(x)*\frac{1}{\sqrt{1-(\frac{x}{t})^2}}.
\end{equation*}
Given that the value of $\rho(x,t)$ is less than zero in this scenario, we exclude this case from our consideration.

In conclusion, we have
\begin{equation}\label{onedimensionrho}
    \rho(x,t) = 2\rho_{in}*\frac{1}{t}F(\frac{x}{t}) = \frac{1}{\pi t}\rho_{in}(x)*\frac{1}{\sqrt{1-(\frac{x}{t})^2}}.
\end{equation}

\begin{figure}
    \centering
    \includegraphics[width=0.7\linewidth]{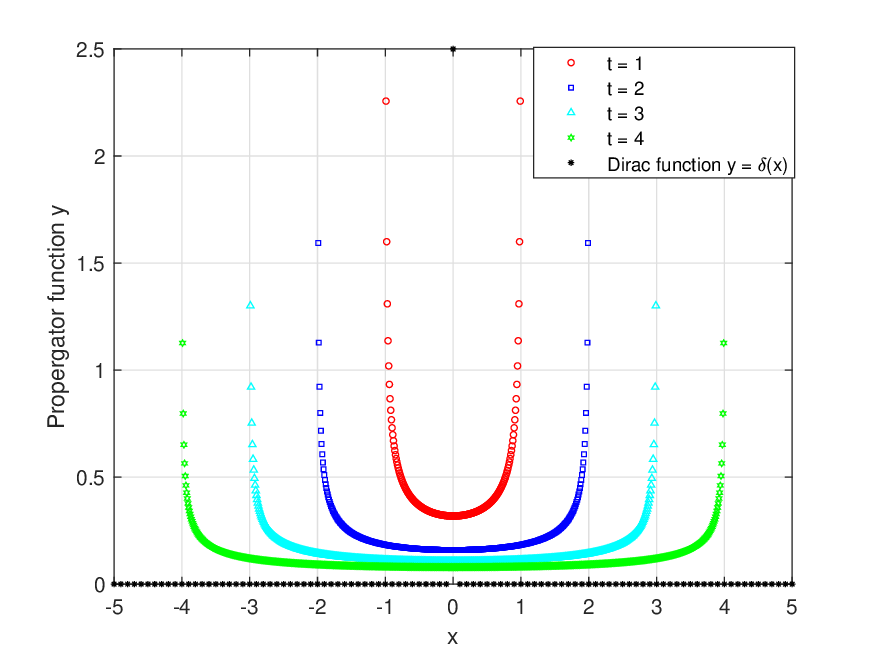}
    \caption{$\rho(x,t)$ defined in \eqref{onedimensionrho} with the initial value $\rho_{in} = \delta(x)$.}
    \label{onedim}
\end{figure}

\paragraph{The numerical results}
To verify our analytical result \eqref{onedimensionrho}, we set the initial value $\rho_{in}(x) = \delta(x)$. Then, $\rho(x,t)$ becomes 
\begin{equation}
    \rho(x,t) = \rho_{in}*\frac{1}{\pi t}F(\frac{x}{t}) = \frac{1}{\pi t}\frac{1}{\sqrt{1-(\frac{x}{t})^2}}.
\end{equation}
For the MSD of this one-dimensional case, we have
\begin{equation}
    \langle x^2(t)\rangle = \int_{\mathbb{R}} x^2 \rho(x,t) \, dx =\frac{t^2}{\pi} \int_{-1}^{1} \frac{y^2}{\sqrt{1-y^2}} \, dy = \frac{1}{2}t^2,
\end{equation}
which signifies the ballistic transport of the particle swarm.

we plotted the theoretical density function $\rho(x,t)$ for $t = 1, 2, 3, 4$. The numerical simulation result is shown in Figure \ref{onedim}. We can observe that the ballistic fronts propagate to both the left and right sides at a speed of $v = \pm 1$. The densities have integrable divergences at the ballistic fronts due to the conservation of the total number of particles. Moreover, the initial condition is gradually forgotten over time, and the solution approaches the universal self-similar profile of the corresponding Green's function.

\subsection{The derivation of the diffusion limit equation of $\rho$ when $\gamma \leq 0$.}\label{derivationrho}
In this subsection, we prove the main theorem (b) for $\gamma \leq 0$. When $\gamma \leq 0$, we set $\mu = 1$. Integrating \eqref{f+rhoin} and \eqref{f-rhoin}, we get 
\begin{equation}\label{rhogamma<=0}
    \Tilde{\hat{\rho}}_{\epsilon} = \int_{\mathbb{V}} \int_{\mathbb{R}} \Tilde{\hat{f}}_{\epsilon} \, dm \, dv =\int_{\mathbb{V}} \left(\int_0^{+\infty} \Tilde{\hat{f}}_{\epsilon,+} \, dm + \int_{-\infty}^0 \Tilde{\hat{f}}_{\epsilon,-} \, dm\right) \, dv = \frac{\epsilon^2 \hat{\rho}_{in}  \int_{\mathbb{V}} \frac{1}{-\lambda(\nu-\lambda -2\epsilon^{\gamma})} + \frac{1}{\nu(\nu-\lambda -2\epsilon^{\gamma})} \, dv}{1-\int_{\mathbb{V}} \frac{\, dv}{-\lambda(\nu-\lambda -2\epsilon^{\gamma})}}.
\end{equation}

\subsubsection{The case $\gamma = 0$}
When $\gamma = 0$, the asymptotic expansions for $\lambda$ and $\nu$ are 
\begin{equation*}
    \lambda = -\frac{\sqrt{5}+1}{2} - \epsilon^2 \left(\frac{1}{\sqrt{5}} s + \frac{1}{5\sqrt{5}} |\xi \cdot v|^2\right) + O(\epsilon^4)-i\left[\frac{1}{\sqrt{5}} \epsilon \xi \cdot v + O(\epsilon^3)\right].
\end{equation*}
\begin{equation*}
    \nu = 1+ \epsilon^2(s + |\xi \cdot v|^2) +  O(\epsilon^4) + i\left[\epsilon \xi \cdot v +O(\epsilon^3)\right].
\end{equation*}

\paragraph{The asymptotic expansion of $\int_{\mathbb{V}} \frac{\, dv}{-\lambda(\nu-\lambda-2)} + \int_{\mathbb{V}} \frac{\, dv}{\nu(\nu-\lambda -2)}$.}
From the asymptotic expansions of $\lambda$ and $\nu$, we have
\begin{equation*}
    -\lambda(\nu-\lambda-2) = 1+ \epsilon^2 \left(\frac{3+\sqrt{5}}{2} s + \frac{5+3\sqrt{5}}{10} |\xi \cdot v|^2\right) + O(\epsilon^4) + i \left[\epsilon \frac{3+\sqrt{5}}{2}  \xi \cdot v + O(\epsilon^3)\right].
\end{equation*}

\begin{equation*}
    \nu(\nu-\lambda-2) = \frac{\sqrt{5}-1}{2} + \epsilon^2 \left(\frac{5+7\sqrt{5}}{10} s + \frac{17\sqrt{5}-25}{50} |\xi \cdot v|^2\right) + O(\epsilon^4) + i \left[\epsilon \frac{5+7\sqrt{5}}{10} \xi \cdot v + O(\epsilon^3)\right].
\end{equation*}
It follows that 
\begin{equation}\label{noint1/lambdagamma0}
    \frac{1}{-\lambda(\nu-\lambda-2)} = 1 -\epsilon^2 \left(\frac{3+\sqrt{5}}{2} s + \frac{20+9\sqrt{5}}{5} |\xi \cdot v|^2 \right) + O(\epsilon^4) - i \left[\epsilon \frac{3+\sqrt{5}}{2}  \xi \cdot v + O(\epsilon^3)\right].
\end{equation}
and
\begin{equation}\label{noint1/nugamma0}
    \frac{1}{\nu(\nu-\lambda-2)} = \frac{1+\sqrt{5}}{2} - \epsilon^2 \left(\frac{25+13\sqrt{5}}{10} s + \frac{405+218\sqrt{5}}{50} |\xi \cdot v|^2 \right) + O(\epsilon^4) - i \left[ \epsilon \frac{25+13\sqrt{5}}{10}  \xi \cdot v + O(\epsilon^3)\right].
\end{equation}

Using \eqref{noint1/lambdagamma0} and \eqref{noint1/nugamma0}, we obtain that
\begin{equation*}
    \int_{\mathbb{V}} \frac{1}{-\lambda(\nu-\lambda-2)} \, dv = 1 - \epsilon^2 \left(\frac{3+\sqrt{5}}{2} s - \frac{20+9\sqrt{5}}{5} \int_{\mathbb{V}}|\xi \cdot v|^2 \, dv \right)+ O(\epsilon^4) + iO(\epsilon^3),
\end{equation*}
and
\begin{equation*}
    \int_{\mathbb{V}} \frac{1}{\nu(\nu-\lambda-2)} \, dv = \frac{1+\sqrt{5}}{2} -\epsilon^2 \left(\frac{25+13\sqrt{5}}{10} s + \frac{405+218\sqrt{5}}{50} \int_{\mathbb{V}}|\xi \cdot v|^2 \, dv \right)+ O(\epsilon^4) + iO(\epsilon^3).
\end{equation*}
As a result, we have
\begin{align}\label{gamma=0fenzi}
    \int_{\mathbb{V}} \frac{\, dv}{-\lambda(\nu-\lambda -2)} + \int_{\mathbb{V}} \frac{\, dv}{\nu(\nu-\lambda -2)} &={} \frac{3+\sqrt{5}}{2} - \epsilon^2 \left(\frac{20+9\sqrt{5}}{5}s + \frac{205+128\sqrt{5}}{50} \int_{\mathbb{V}} |\xi \cdot v|^2 \, dv\right) \notag \\
    &+{} O(\epsilon^4) + iO(\epsilon^3).
\end{align}

\paragraph{The asymptotic expansion of $\left[1-\int_{\mathbb{V}} \frac{\, dv}{-\lambda(\nu-\lambda -2)}\right]^{-1}$.}
By calculation, 
\begin{equation}\label{gamma=0fenmu}
    1-\int_{\mathbb{V}} \frac{\, dv}{-\lambda(\nu-\lambda -2)} = \epsilon^2 \left(\frac{3+\sqrt{5}}{2}s + \frac{20+9\sqrt{5}}{5} \int_{\mathbb{V}} |\xi \cdot v|^2 \, dv\right) + O(\epsilon^4) + i O(\epsilon^3).
\end{equation}
Then,
\begin{equation*}
    \left[1-\int_{\mathbb{V}} \frac{\, dv}{-\lambda(\nu-\lambda -2)} \right]^{-1} = \epsilon^{-2}\left(\frac{3+\sqrt{5}}{2}s + \frac{20+9\sqrt{5}}{5} \int_{\mathbb{V}} |\xi \cdot v|^2 \, dv\right)^{-1} + O(1) + i O(\epsilon^{-1}).
\end{equation*}
Substituting \eqref{gamma=0fenzi} and \eqref{gamma=0fenmu} into \eqref{rhogamma<=0}, one can get
\begin{align*}
    \Tilde{\hat{\rho}}_{\epsilon} &={} \hat{\rho}_{in}\left[\frac{3+\sqrt{5}}{2}+O(\epsilon^2) + iO(\epsilon^3)\right]\left[\left(\frac{3+\sqrt{5}}{2}s + \frac{20+9\sqrt{5}}{5} \int_{\mathbb{V}} |\xi \cdot v|^2 \, dv\right)^{-1} + O(\epsilon^2) + iO(\epsilon)\right] \notag \\
    &={} \hat{\rho}_{in} \left[\frac{3+\sqrt{5}}{2} \left(\frac{3+\sqrt{5}}{2}s + \frac{20+9\sqrt{5}}{5} \int_{\mathbb{V}} |\xi \cdot v|^2 \, dv\right)^{-1} + O(\epsilon^2) + iO(\epsilon) \right].
\end{align*}
Letting $\epsilon \to 0$, we get the limit equation of $\Tilde{\hat{\rho}} = \lim_{\epsilon \to 0} \Tilde{\hat{\rho}}_{\epsilon}$ writing
\begin{equation*}
    \Tilde{\hat{\rho}} = \hat{\rho}_{in} \left[\frac{3+\sqrt{5}}{2} \left(\frac{3+\sqrt{5}}{2}s + \frac{20+9\sqrt{5}}{5} \int_{\mathbb{V}} |\xi \cdot v|^2 \, dv\right)^{-1}\right].
\end{equation*}
Multiply both sides by $\frac{2}{3+\sqrt{5}} \left(\frac{3+\sqrt{5}}{2}s + \frac{20+9\sqrt{5}}{5} \int_{\mathbb{V}} |\xi \cdot v|^2 \, dv\right)$,
one gets
\begin{equation*}
    s \Tilde{\hat{\rho}} - \hat{\rho}_{in} + \Tilde{\hat{\rho}} |\xi|^2 \frac{15+7\sqrt{5}}{10} \int_{\mathbb{V}} (e_{\xi} \cdot v)^2 \, dv = 0,
\end{equation*}
where $\xi = |\xi| e_{\xi}$. After the inverse Laplace and Fourier transform, we get the limit equation of $\rho(x,t)$
\begin{equation*}
    \partial_t \rho(x,t) - C\Delta_x \rho(x,t) = 0,
\end{equation*}
where the constant $C$ depends on the dimension $d$. By calculation, for $d \geq 2$, $C = \frac{15+7\sqrt{5}}{10} \int_{\mathbb{V}} (e_{\xi} \cdot v)^2 \, dv = \frac{15+7\sqrt{5}}{30}$, and for $d = 1$, $C = \frac{15+7\sqrt{5}}{10}$.
\subsubsection{The case $\gamma< 0$}
When $-1 \leq \gamma<0$, the expansion of $\lambda$ is given by
\begin{equation}\label{lambdacase1}
     \lambda =-\epsilon^{\gamma}-\epsilon^{-\gamma} - O(\epsilon^{-3\gamma}) - i\left[\epsilon^{1-\gamma} \xi \cdot v + O(\epsilon^{1-2\gamma})\right].
\end{equation}
When $\gamma < -1$, the expansion is 
\begin{equation*}
    \lambda = -\epsilon^{\gamma}-\epsilon^{-\gamma} - O(\epsilon^{2-\gamma})-i\left[\epsilon^{1-\gamma} \xi \cdot v + O(\epsilon^{1-2\gamma})\right].
\end{equation*}
Both scenarios yield consistent results, so we focus on the first case \eqref{lambdacase1}. As for $\nu$, the expansion is given by
\begin{equation*}
    \nu = \epsilon^{\gamma}+ \epsilon^{2-\gamma}s +  O(\epsilon^{2-3\gamma}) + i\left[\epsilon^{1-\gamma} \xi \cdot v +O(\epsilon^{3-2\gamma})\right].
\end{equation*}
Next, we derive the asymptotic expansion of \eqref{rhogamma<=0}.

\paragraph{The asymptotic expansion of $\int_{\mathbb{V}} \frac{\, dv}{-\lambda(\nu-\lambda-2\epsilon^{\gamma})} + \int_{\mathbb{V}} \frac{\, dv}{\nu(\nu-\lambda -2\epsilon^{\gamma})}$.}
It can be obtained that
\begin{equation*}
    \nu-\lambda-2\epsilon^{\gamma} = \epsilon^{-\gamma} + O(\epsilon^{-3\gamma}) + i\left[\epsilon^{1-\gamma} 2\xi \cdot v + O(\epsilon^{1-2\gamma})\right].
\end{equation*}
Then, we have
\begin{equation*}
    -\lambda(\nu-\lambda-2\epsilon^{\gamma}) = 1+ O(\epsilon^{-2\gamma}) + i\left[\epsilon 2\xi \cdot v + O(\epsilon^{1-\gamma})\right].
\end{equation*}
\begin{equation*}
    \nu(\nu-\lambda-2\epsilon^{\gamma}) = 1 + O(\epsilon^{-2\gamma}) + i\left[\epsilon 2\xi \cdot v + O(\epsilon^{1-\gamma})\right].
\end{equation*}
Thus the inverses are
\begin{equation*}
    \frac{1}{-\lambda(\nu-\lambda-2\epsilon^{\gamma})} = 1 + O(\epsilon^{-2\gamma}) + O(\epsilon^2) - i\left[\epsilon 2\xi \cdot v + O(\epsilon^{1-\gamma})\right]. 
\end{equation*}
\begin{equation*}
    \frac{1}{\nu(\nu-\lambda-2\epsilon^{\gamma})} = 1 + O(\epsilon^{-2\gamma}) - i\left[\epsilon 2\xi \cdot v + O(\epsilon^{1-\gamma})\right].
\end{equation*}

Consequently, it can be obtained 
\begin{equation}\label{fenzigammageq-1}
    \int_{\mathbb{V}} \frac{\, dv}{-\lambda(\nu-\lambda-2\epsilon^{\gamma})} + \int_{\mathbb{V}} \frac{\, dv}{\nu(\nu-\lambda -2\epsilon^{\gamma})} = 2 + O(\epsilon^{-2\gamma}) + iO(\epsilon^{1-\gamma}).
\end{equation}

\paragraph{The asymptotic expansion of $1-\int_{\mathbb{V}} \frac{\, dv}{-\lambda(\nu-\lambda-2\epsilon^{\gamma})}$.}

From \eqref{eigeqnlambda} and \eqref{eigeqnnu}, we derive
\begin{equation*}
    \lambda^2 = -\epsilon^{\gamma} \lambda + \epsilon^2 s+ i \epsilon \xi \cdot v + 1,
\end{equation*}
\begin{equation*}
    \nu^2 = \epsilon^{\gamma} \nu + \epsilon^2 s+ i \epsilon \xi \cdot v.
\end{equation*}
As a result, for $\frac{1}{-\lambda(\nu-\lambda-2\epsilon^{\gamma})}$, we can get
\begin{align*}
    \frac{1}{-\lambda(\nu-\lambda-2\epsilon^{\gamma})} &={} \frac{\nu}{-\lambda\nu(\nu-\lambda-2\epsilon^{\gamma})} = \frac{\nu}{-\lambda \nu^2 + \lambda^2 \nu + 2\epsilon^{\gamma} \lambda \nu} \notag \\
    &={} \frac{\nu}{-\lambda(\epsilon^{\gamma} \nu + \epsilon^2 s+ i \epsilon \xi \cdot v) + \nu(-\epsilon^{\gamma} \lambda + \epsilon^2 s+ i \epsilon \xi \cdot v + 1) + 2\epsilon^{\gamma} \lambda \nu} \notag \\
    &={} \frac{\nu}{\nu + (\nu-\lambda)(\epsilon^2 s+ i \epsilon \xi \cdot v)} = \frac{1}{1+\frac{\nu-\lambda}{\nu}(\epsilon^2 s+ i \epsilon \xi \cdot v)}.
\end{align*}

By calculation, we obtain
\begin{equation*}
    \nu -\lambda = 2\epsilon^{\gamma} + \epsilon^{-\gamma} + O(\epsilon^{-3\gamma}) + i\left[\epsilon^{1-\gamma} 2 \xi \cdot v + O(\epsilon^{1-2\gamma})\right].
\end{equation*}
\begin{equation*}
    \frac{1}{\nu} = \epsilon^{-\gamma} - \epsilon^{2-3\gamma}s + O(\epsilon^{2-5\gamma}) - i\left[\epsilon^{1-3\gamma} \xi \cdot v + O(\epsilon^{1-4\gamma})\right].
\end{equation*}
\begin{equation*}
    \frac{\nu-\lambda}{\nu} = 2+\epsilon^{-2\gamma} +O(\epsilon^{-4\gamma}) - iO(\epsilon^{1-3\gamma}).
\end{equation*}
\begin{equation*}
    1+\frac{\nu-\lambda}{\nu}(\epsilon^2 s + i\epsilon \xi \cdot v) = 1+ \epsilon^2 2s + O(\epsilon^{2-2\gamma}) + i\left[\epsilon 2 \xi \cdot v + O(\epsilon^{1-2\gamma})\right].
\end{equation*}
Then,
\begin{align*}
    \frac{1}{-\lambda(\nu-\lambda-2\epsilon^{\gamma})} &={} \left[ 1+\frac{\nu-\lambda}{\nu}(\epsilon^2 s + i\epsilon \xi \cdot v)\right]^{-1} \notag \\
    &={} 1-\epsilon^2 (2s + 4|\xi \cdot v|^2) + O(\epsilon^{2-2\gamma}) - i\left[\epsilon 2\xi \cdot v + O(\epsilon^{1-2\gamma})\right].
\end{align*}
As a result, the integral becomes
\begin{equation*}
    \int_{\mathbb{V}} \frac{1}{-\lambda(\nu-\lambda-2\epsilon^{\gamma})} \, dv = 1- \epsilon^2 \left(2s + \int_{\mathbb{V}} 4 |\xi \cdot v|^2 \, dv \right) + O(\epsilon^{2-2\gamma}).
\end{equation*}
Hence, it can be derived that
\begin{equation}\label{fenmugammageq-1}
    1-\int_{\mathbb{V}} \frac{1}{-\lambda(\nu-\lambda-2\epsilon^{\gamma})} \, dv = \epsilon^2 \left(2s + \int_{\mathbb{V}} 4 |\xi \cdot v|^2 \, dv \right) + O(\epsilon^{2-2\gamma}).
\end{equation}

\paragraph{The derivation of the limit equation.}
From \eqref{derivationrho}, it is straightforward to infer
\begin{equation}\label{newderivationrho}
    \Tilde{\hat{\rho}}_{\epsilon} \left[1-\int_{\mathbb{V}} \frac{1}{-\lambda(\nu-\lambda-2\epsilon^{\gamma})} \, dv\right] = \epsilon^2 \hat{\rho}_{in} \left[\int_{\mathbb{V}} \frac{\, dv}{-\lambda(\nu-\lambda-2\epsilon^{\gamma})} + \int_{\mathbb{V}} \frac{\, dv}{\nu(\nu-\lambda -2\epsilon^{\gamma})}\right].
\end{equation}
Substituting \eqref{fenzigammageq-1} and \eqref{fenmugammageq-1} into \eqref{newderivationrho} and dividing both sides by $2\epsilon^2$ simultaneously, one has
\begin{equation}\label{formaleqngamma<0}
    \Tilde{\hat{\rho}}_{\epsilon} \left(s + \int_{\mathbb{V}} 2 |\xi \cdot v|^2 \, dv \right) = \hat{\rho}_{in}(1+O(\epsilon^{-2\gamma})).
\end{equation}
Letting $\epsilon \to 0$ simultaneously on both sides of \eqref{formaleqngamma<0}, we obtain 
\begin{equation*}
    s \Tilde{\hat{\rho}} - \hat{\rho}_{in} + \Tilde{\hat{\rho}} |\xi|^2 2 \int_{\mathbb{V}} e_{\xi} \cdot v \, dv,
\end{equation*}
where $\xi = |\xi| e_{\xi}$. After the inverse Laplace and Fourier transform, the limit equation of $\rho(x,t)$ is
\begin{equation*}
 \partial_t \rho(x,t) - C \Delta_x \rho(x,t) = 0,
\end{equation*}
where $C$ depends on dimension $d$. By computation, for dimension $d \geq 2$, $C = 2\int_{\mathbb{V}} (e_{\xi} \cdot v)^2 \, dv = \frac{2}{3}$, and for $d = 1$, $C = 2$.

\section{Conclusion and discussions}
In this paper, we have investigated the crossover phenomenon from ballistic transport to normal diffusion induced by intracellular adaptation and noise. Through a multi-scale analysis of a kinetic model with internal state, we have successfully captured and explained this transition behavior. Our findings highlight the crucial role of adaptation time in modulating the bacterial search strategy, with longer adaptation times favoring ballistic motion for efficient exploration and shorter adaptation times leading to normal diffusion for localized exploitation.

The mathematical analysis and numerical simulations presented herein provide a robust framework for understanding the complex interplay between intracellular dynamics and macroscopic movement patterns. The agreement between the individual based model and kinetic model further strengthens the validity of our conclusions. The construction of this model is also partially inspired by recent studies of the firing mechanism of neurons \cites{zhou2021investigating, liu2022rigorous}. While our models have provided valuable insights, it is important to acknowledge their limitations. We have employed simplified representations of intracellular signaling and assumed a constant velocity for bacterial runs. Future studies could incorporate more detailed biochemical pathways and consider variations in run speeds. Additionally, exploring the impact of external stimuli, such as chemical attractant gradients, on the crossover phenomenon would be a valuable extension of this work.


\bigskip
\textbf{Acknowledgement:} 
Z. Xue and M. Tang  are
 supported by NSFC12031013 and the Strategic Priority Research Program of Chinese Academy of Sciences Grant No.XDA25010401; Z. Zhou is supported by the National Key R\&D Program of China with project number 2021YFA1001200, and the NSFC with grant number 12171013.

\appendix
\section{The numerical simulation in Subsection \ref{two-state}}\label{twostateappendix}
The $n$th sample is represented by the CheY-P concentration $Y^n$ and the tumbling state $s^n$. When $s^n = 1$, the sample rotates counterclockwise (CCW). When $s^n = -1$, the sample is in the clockwise (CW) state. In our simulation, $N = 100$ samples are tracked, and each sample evolves according to the following algorithm with a time step $\Delta \tau = 0.1$. $\Lambda_1$ and $\Lambda_2$ are defined in Equation \eqref{TupaperLambda}. We denote the numerical approximations of $Y^n(k \Delta \tau)$, $s^n(k \Delta \tau)$, $\Lambda_1^n(k \Delta \tau)$, and $\Lambda_2^n(k \Delta \tau)$ by $Y_k^n$, $s_k^n$, $\Lambda_{1,k}^n$, and $\Lambda_{2,k}^n$, respectively.

\paragraph{Initialization.} We set the correlation time $T_m = 6000$, the mean value of the CheY-P concentration $\bar{Y} = 5$, and $\sigma = 0.456$. The tumbling rate is given by Equation \eqref{TupaperLambda}. We set $\alpha_1 = 10$ and $\alpha_2 = -2$. We assign the characteristic switching times $t_0 = 300$ and $t_1 = 30$. For $n = 1, \ldots, N$, the initial value of $Y^n$ is $\bar{Y}$, and the initial state of a sample is set randomly to either CW or CCW, which means $s^n$ is randomly set to $1$ or $-1$. We evolve $(Y^n, \Lambda_1^n, \Lambda_2^n)$ as follows:

\paragraph{Time evolution.} For each time step $k$ (=1 initially), we perform the following calculations repeatedly until $k = 6 \times 10^6$:
\begin{itemize}
    \item[1)] Update $Y_k^n$. Generate a random number denoted by $\Delta B_t$ using the normal distribution $N(0,\Delta \tau)$. Use the Euler-Maruyama method to update $Y^n$
    \begin{equation*}
        Y_{k+1}^n = Y_k^n - \frac{Y_k^n - \bar Y}{T_m} \Delta \tau + \sigma \Delta B_t.
    \end{equation*}
    \item[2)] Update the tumble rate
    \begin{equation*}
        \Lambda_1^n = t_0^{-1} \exp(-\alpha_1 \frac{Y_k^n - \bar Y}{\bar Y})
    \end{equation*}
    \begin{equation*}
        \Lambda_2^n = t_1^{-1} \exp(-\alpha_2 \frac{Y_k^n - \bar Y}{\bar Y})
    \end{equation*}
    \item[3)] If $s^n=1$, set the jumping rate $\Lambda = \Lambda_1$.
    \par
    Otherwise, $s^n = -1$, set $\Lambda = \Lambda_2$.
    \item[4)] Generate a random number $r^n$ uniformly distributed in $[0,1]$,
    \par
    If $r^n \leq \Lambda \Delta \tau \exp(\Lambda \Delta \tau)$, set $s^n = -s^n$,
    \par
    Otherwise, set $s^n = s^n$.   
    
    \item[5)]Otherwise, set $k = k+1$, and go back to Step 1).
\end{itemize}
\section{The numerical simulation of the one-state model}\label{one-statemodelappendix}
The $n$th sample is represented by the internal state $m^n$ and the integral of $\Lambda$, denoted by $I^n$. We denote the numerical approximations of $m^n(k \Delta \tau)$ and the integral $I^n(k \Delta \tau)$ by $m_k^n$ and $I_k^n$, respectively. In our simulation, $N = 100$ samples are tracked, and each sample evolves according to the following algorithm with a time step $\Delta \tau = 0.1$.

\paragraph{Initialization.} The initial values $m^n$ and $I^n$ for all samples are $0$. For each sample, the total running time $T_t$ is $6\times10^5$. The value of $T_m = 6000$ and $\sigma = 0.456$. $\Lambda$ and $g$ are, respectively, in \eqref{Lambda} and \eqref{gm}. We generate a series of independent random numbers $\Gamma \sim \exp(1)$. For $n = 1,...,N$, we evolve $(m^n,I^n)$ as follows:

\paragraph{Time evolution.} For each time step $k$ (=1 initially), we perform the following calculations repeatedly until $k=6 \times 10^6$:
\begin{itemize}
    \item[1)] Update the internal state. Generate a random number denoted by $\Delta B_t$ using the normal distribution $N(0,\Delta \tau)$. Use the Euler-Maruyama method to update $X^n$
    \begin{equation*}
        m_{k+1}^n = m_k^n - \frac{1}{T_m} g(m_k^n) \Delta \tau + \sqrt{2} \Delta B_t.
    \end{equation*}
    \item[2)] Update $I^n$. Update the value of $I^n$ by
    \begin{equation*}
        I_{k+1}^n = I_k^n + \Lambda(m_{k+1}^n) \Delta \tau.
    \end{equation*}
    \item[3)] If $I_k^n \geq \Gamma$, set $m_{k+1}^i = 0$ and $I_k^n = 0$. 
    
    \item[4)] Otherwise, set $k \gets k+1$, and go back to Step 1).
\end{itemize}

\section{The algorithm for the kinetic model with cell movement}\label{IBMappendix}
Each sample is represented by $(x^n,m^n,v^n,I^n)$. We denote the numerical approximations of $x^n(k \Delta \tau)$, $m^n(k \Delta \tau)$, $v^n(k \Delta \tau)$ and $I^n(k \Delta \tau)$ by $x_k^n$, $m_k^n$, $v_k^n$ and $I_k^n$, respectively. In our simulation, $N = 1000$ samples are tracked and each sample evolves by the following algorithm with the time step $\Delta \tau= 1$ for $T_m = 1, 10, 100, 1000$ and $+\infty $. When $T_m = 10^{-2}$, $\Delta \tau = 10^{-2}$ Here we show the details of the algorithm.

\paragraph{Initialization.} The initial values $x^n$, $m^n$and $I^n$ for all samples are $0$. We fix $V_0 = 0.02$. The initial velocity $v^n$ is randomly set to $V_0$ or $-V_0$ with equal probability. For each sample, the total running time $T_t$ is $10^6$. The value of $T_m$ is set to $10^{-2}$, $1$, $10$, $100$, $1000$ and $+\infty$. $\Lambda$ and $g$ are, respectively, in \eqref{Lambda} and \eqref{gm}. We generate a series of independent random numbers $\Gamma \sim \exp(1)$. For $n = 1,...,N$, we evolve $(x^n,m^n,v^n,I^n)$ as follows:

\paragraph{Time evolution.} For each time step $k$ (=1 initially), we perform the following calculations repeatedly until $k=T_t/\Delta \tau$:
\begin{itemize}
    \item[1)] Update the internal state. Generate a random number denoted by $\Delta B_t$ using the normal distribution $N(0,\Delta \tau)$. Use the Euler-Maruyama method to update $X^n$
    \begin{equation*}
        m_{k+1}^n = m_k^n - \frac{1}{T_m} g(m_k^n) \Delta \tau + \sqrt{2} \Delta B_t.
    \end{equation*}
    \item[2)] Update $I^n$. Update the value of $I^n$ by
    \begin{equation*}
        I_{k+1}^n = I_k^n + \Lambda(m_{k+1}^n) \Delta \tau.
    \end{equation*}
    \item[3)] Update $x^n$. Set $x^n_{k+1} = x^n_{k} + v_k^n \Delta \tau$.
    \item[4)] If $I_k^n \geq \Gamma$, set $m_{k+1}^i = 0$ and $I_k^n = 0$. Generate a random number $r$ uniformly distributed in $[0,1]$. 
     
    \qquad If $r \leq \frac{1}{2}$, set $v_{k+1}^i = -v_k^i$. Otherwise, set $v_{k+1}^i = v_k^i$.
    
    \item[5)] Otherwise, set $k \gets k+1$, and go back to Step 1).
\end{itemize}

\bibliographystyle{unsrt}
\bibliography{ref.bib}

\end{document}